\documentclass[a4paper,pra,aps,12pt,showkeys,nofootinbib]{revtex4}

\usepackage[T1]{fontenc}
\usepackage[utf8]{inputenc}
\usepackage{units}
\usepackage{textcomp}
\usepackage{amsmath}    
\usepackage{graphicx}   
\usepackage{verbatim}   
\usepackage{color}      
\usepackage{subfigure}  
\usepackage{hyperref}   
\hypersetup{citebordercolor=.2 .8 .2,urlbordercolor=0 .75 1}
\usepackage{nicefrac}
\usepackage{amsfonts}
\usepackage{amsthm}
\usepackage{amssymb}
\usepackage{booktabs}
\usepackage{fancyhdr}
\usepackage{moreverb}
\usepackage[brazil]{babel}

\makeatletter
\renewcommand*{\email}[1][E-mail: ]{\begingroup\sanitize@url\@email{#1}}
\makeatother

\DeclareRobustCommand{\quatrorevertido}{\reflectbox{$ \boldsymbol{4} $}}
\DeclareRobustCommand{\cincorevertido}{\reflectbox{$ \boldsymbol{5} $}}
\DeclareRobustCommand{\noverevertido}{\reflectbox{$ \boldsymbol{9} $}}
\DeclareRobustCommand{\seterevertido}{\reflectbox{$ \boldsymbol{7} $}}
\DeclareRobustCommand{\doisrevertido}{\reflectbox{$ \boldsymbol{2} $}}

\begin{document}

	\title{Uma história e o uso de espelhos no aprendizado das\\regras de sinais envolvendo números negativos}

	\author{M. F. Araujo de Resende}
	\email{resende@if.usp.br}
	
	\affiliation{Instituto de Física, Universidade de São Paulo, 05508-090 São Paulo SP, Brasil}

	\date{\today }
	
	\begin{abstract}
		A proposta deste artigo (que, em muitas passagens, eu tomo a liberdade de escrever na primeira pessoa do singular por trazer uma experiência pessoal) é a de mostrar como os professores de Física e de Matemática do Ensino Básico podem ensinar, através de uma atividade simples e lúdica, as regras de sinais que estão relacionadas para com as operações aritméticas básicas envolvendo números positivos e negativos. Como? Usando uma régua, ou qualquer outro objeto que possa ser interpretado como uma reta numérica, e um espelho plano. Aliás, como a justificativa do porquê esse método funciona é física e matemática, este artigo traz uma apresentação bastante sucinta sobre a formação de imagens por espelhos planos e esféricos posto que alguns professores de Matemática podem (apenas como uma possibilidade) não saber como isso acontece. Por fim, diante na natureza interdisciplinar desta atividade, algumas considerações adicionais também são apresentadas sobre como esses professores podem usá-la para instigar os estudantes a observarem como a Física e a Matemática, juntas, fazem parte do dia a dia de todos nós.
	\end{abstract}



	

	\maketitle

	\section{\label{intro}Era uma vez numa sala de aula...}
	
		Se existe uma coisa que pode ser bem difícil para alguns estudantes é entender o porquê daquelas famigeradas regras de sinais que estão relacionadas às operações básicas envolvendo números positivos e negativos.
		\begin{quote}
			\emph{\textquotedblleft Mais com mais é mais,\\ Mais com menos é menos,\\ Menos com mais também é menos e\\ Menos com menos é mais.\\ Por que isso é assim?\textquotedblright}
		\end{quote}
		Ao menos era essa a pergunta que não saía da minha cabeça aos doze anos de idade, algum tempo depois da minha introdução aos números negativos.
		
		Confesso que eu não posso dizer que eu tive \textquotedblleft grandes\textquotedblright \hspace*{0.01cm} dificuldades com a Matemática ao longo da minha vida já que, desde muito cedo, as minhas notas eram muito altas e eu já gostava muito dela. E como esse meu gosto era tanto, eu chego a desconfiar que algumas pessoas, que me conheceram na época em que eu conheci os primeiros números e as suas primeiras somas, deviam achar que eu era uma criança bem estranha em algum sentido. Afinal de contas, eu sempre arranjava um jeito de me cercar de números por todos os lados a ponto de, por exemplo, eu sempre carregar, comigo, uma caneta e um caderno (inclusive durante as férias) só para eu ficar escrevendo números e inventando algumas somas com as letras do alfabeto.
		
		Aliás, ver números e letras sendo somados era muito lindo para mim, e a certeza que eu já tinha, desde os meus seis anos de idade, era a de que, quando eu \textquotedblleft crescesse\textquotedblright , eu iria saber fazer todas as contas possíveis, principalmente todas aquelas que permitiam entender as coisas do Universo. Só que, apesar de eu realmente não poder afirmar que eu tive \textquotedblleft grandes\textquotedblright \hspace*{0.01cm} dificuldades com a Matemática, eu já não posso dizer o mesmo em relação à maioria dos professores que me deram aulas de Matemática. A Matemática que eles me apresentavam era muito chata já que, nas suas aulas, eles nunca propunham qualquer desafio interessante, eles nunca tentavam construir uma correspondência daquela Matemática para com a realidade. Basicamente a única coisa que eles faziam era pedir para decorar tabuadas e um monte de regras para fazer operações sem nunca justificá-las. E como aquela Matemática não era a que eu já gostava, com o passar do tempo eu fui desanimando, eu fui deixando de fazer todas as tarefas que não valiam notas e eu só fazia essas, que valiam notas, no \textquotedblleft último minuto\textquotedblright \hspace*{0.01cm} para minimizar toda aquela chatice na minha vida. O meu ânimo só voltou, de verdade, quando eu estava com os mesmos doze anos de idade de quando eu conheci os números negativos, principalmente depois que todas as aulas começaram a ser voltadas para as operações envolvendo letras e números.
			
		Embora explicar toda a lógica que existe por trás da soma de todas as frações seja algo até que bem simples de ser feito por um professor, esse mesmo cenário parece não se aplicar quando, por exemplo, esse mesmo professor precisa explicar por que as multiplicações precisam ser efetuadas antes das somas, ou mesmo por que \textquotedblleft todo\textquotedblright \hspace*{0.01cm} número elevado a zero é sempre igual a um\footnote{Eu coloquei aspas em \textquotedblleft todo\textquotedblright \hspace*{0.01cm} pois, por mais que a maioria dos professores que me deram aula no Ensino Básico tenham sempre me afirmado que \textquotedblleft zero elevado a zero é sempre zero\textquotedblright \hspace*{0.01cm} (inclusive me dando errado quando eu não associava qualquer resultado a essa operação naquela época), a verdade é que o resultado de zero elevado a zero \textbf{não é bem definido} \cite{elon1,elon2}: afinal de contas, embora alguns algebristas até lidem com a definição de que $ \boldsymbol{0^{0}} = \boldsymbol{1} $, pelo ponto de vista da Análise essa definição não faz muito sentido. Ou seja, eu realmente tinha alguma razão em não associar qualquer resultado a essa operação.}. A maioria dos professores de Matemática do Ensino Básico apenas pede para que os estudantes \textbf{decorem} todas essas regras e, consequentemente, acabam atribuindo notas a esses estudantes apenas pela capacidade que eles têm de cumprí-las e não de entendê-las. Ou seja, pelo ponto de vista pedagógico, é possível afirmar que essa situação pode ser interpretada como um belo exemplo da educação \textquotedblleft bancária\textquotedblright \hspace*{0.01cm} que Paulo Freire menciona na Ref. \cite{bancaria}, onde os professores tratam os estudantes como meros recipientes que sempre precisam ser preenchidos com alguma informação sem dar qualquer espaço para qualquer reflexão.
		
		No meu caso, por exemplo, foram poucos os professores de Matemática que fugiram um pouco dessa educação \textquotedblleft bancária\textquotedblright \hspace*{0.01cm} e conseguiram me explicar os porquês de algumas operações: entre os treze professores que me deram aulas de Matemática no Ensino Básico, apenas \textbf{um} único deles nunca se esquivou de discutir esses porquês comigo. A única coisa que os demais faziam era passar listas imensas, com exercícios bem chatos e repetitivos, onde as palavras de ordem para resolvê-las sempre eram \textquotedblleft decorar\textquotedblright \hspace*{0.01cm} e \textquotedblleft reproduzir\textquotedblright . E foi exatamente isso que aconteceu quando, aos doze anos de idade, eu cheguei na (antiga) sexta série do Ensino Básico e a professora que me dava aulas de Matemática começou a fazer umas contas com números negativos sempre repetindo
		\begin{quote}
			\emph{\textquotedblleft Mais com mais é mais,\\ Mais com menos é menos,\\ Menos com mais é menos e\\ Menos com menos é mais\textquotedblright ,}
		\end{quote}
		sem nunca explicar o porquê desse \textquotedblleft mantra\textquotedblright \hspace*{0.01cm} para ninguém. Aquilo deu um verdadeiro  \textquotedblleft nó\textquotedblright \hspace*{0.01cm} na minha cabeça! A coisa simplesmente não ia! Consequência? Eu, que sempre fui um dos melhores estudantes de todas as escolas por onde eu passei, pela primeira vez na vida me vi numa situação completamente inusitada: eu fiquei de recuperação em Matemática no segundo bimestre daquele ano (já que, a partir daquele ano, foi instituído um período de recuperação no final de cada bimestre).
			
		Só que, no bimestre seguinte, aconteceu uma coisa diferente: essa professora resolveu substituir a primeira prova por uma lista imensa de exercícios, justamente uma daquelas listas que eu odiava fazer. Acho que tinha mais de duzentas contas que precisavam ser calculadas naquela lista. Só que, para piorar um pouco mais a minha situação, essa lista precisava ser feita em grupo.
		
		Confesso que eu não gostava muito de fazer atividades (que valiam notas) em grupo naquela época. E não era nem por um motivo antissocial: o motivo era o simples fato de não ser incomum alguns dos meus colegas abusarem do fato de eu sempre tirar as melhores notas da turma para não fazerem as suas partes nessas atividades. Ou seja, sempre que eu fazia essas atividades em grupo, não era nada incomum eu me ver numa situação onde eu precisava fazer tudo, por conta própria, só para eu não receber uma nota que eu julgava não merecer.
		
		Só que o caso dessa lista conseguiu ser o pior dos piores. Afinal, apesar da ideia dessa professora até ter sido muito boa (já que a sua proposta provavelmente era a de fazer uma atividade colaborativa, onde os estudantes se ajudassem e aprendessem uns com os outros), não se tratava de um trabalho em grupo onde os estudantes podiam formar os seus grupos por conta própria: todos os grupos foram definidos pela própria professora através de alguns sorteios e o sorteio, que acabou definindo o meu grupo, acabou me colocando justamente ao lado dos estudantes mais problemáticos daquela turma, os quais adoravam praticar \textquotedblleft bullying\textquotedblright \hspace*{0.01cm} com todos, principalmente contra mim. Ou seja, além de todos os problemas que eu já tinha com aqueles benditos números negativos, eu ainda precisava lidar com um \textquotedblleft novo\textquotedblright \hspace*{0.01cm} problema, agora de ordem humana. Resultado: eu não consegui fazer nada com aquele grupo, já que, por definição ou por medo dos assediadores, ninguém ali estava disposto a ser colaborativo para comigo.
		
		Aliás, eu sequer tive acesso àquela lista enquanto eu estive naquele grupo já que apenas uma cópia foi fornecida para cada grupo e, no caso da que foi entregue para esse meu grupo, ela foi entregue justamente nas mãos de um dos assediadores. Eu cheguei até a relatar o causo para a professora, pedindo para que ela me mudasse de grupo, já que as pessoas daquele grupo sequer me davam acesso à lista pra eu resolver. Mas ela me disse que ela não poderia fazer isso pois, do contrário, outros componentes, dos outros grupos, acabariam fazendo o mesmo pedido para ela e seria uma verdadeira confusão. Porém, como a minha situação era meio \textquotedblleft atípica\textquotedblright , ela acabou me fornecendo uma cópia extra da lista e deixou que eu entregasse tudo individualmente caso eu conseguisse.
			
\begin{figure}[!t]
			\centering
			\includegraphics[viewport=360 10 0 144,scale=1.3]{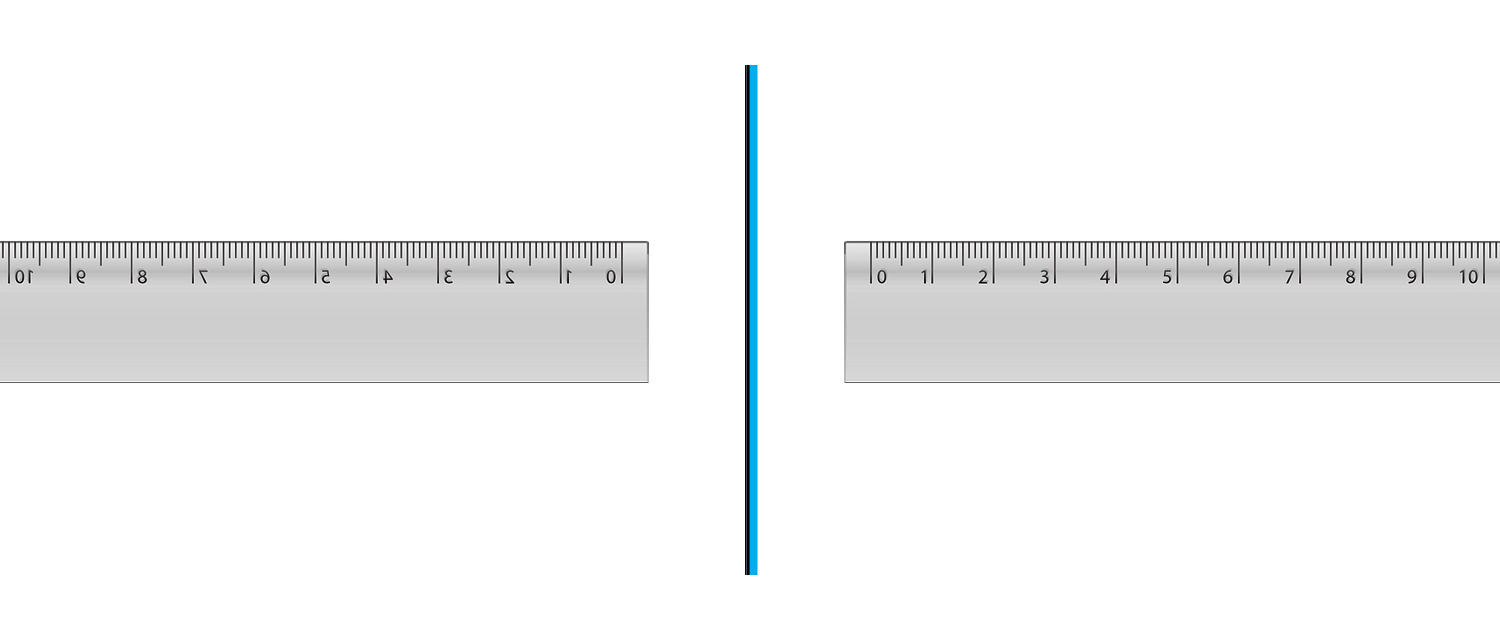}
			\caption{\label{regua-e-imagem}À direita nós podemos ver uma régua que foi colocada diante de um espelho plano (aqui representado pelo segmento de reta destacado em preto), mais especificamente da sua superfície espelhada (que está sendo destacada em ciano). Já à esquerda nós podemos ver a imagem que esse espelho plano foi capaz de construir para essa régua.} 
		\end{figure}

		Eu fiquei quase uma semana sem ir à escola, sozinha, em casa, enquanto a minha mãe estava no trabalho, só tentando resolver aquela bendita lista. Aliás, a minha mãe, que trabalhava na mesma escola onde eu estudava, até comentou num desses dias comigo que, na noite anterior, essa minha professora havia sugerido que a minha mãe contratasse um professor particular de Matemática para ele \textbf{resolver} aquela lista para mim. Ouvindo a minha mãe me dizer aquilo, eu disse
		\begin{quote}
			-- Pra quê?
			
			-- Sabia que você ia falar isso...
			
			-- Pagar alguém pra resolver essa lista, sendo que eu sei que eu sou capaz de fazer tudo sozinha?
			
			-- Foi exatamente isso que eu falei pra ela. \emph{\textquotedblleft Se eu conheço bem a Fernanda, ela jamais vai querer fazer isso. Tudo que ela sempre precisou fazer, ela sempre fez sozinha, ela sempre conseguiu. Deixa ela quieta lá em casa que ela dá um jeito.\textquotedblright }.
			
			-- E pra quê pagar alguém pra resolver essa lista pra mim sendo que sou que tenho que fazer? Sou eu que tenho que aprender! O correto seria ela explicar as coisas direito ou, então, perguntar se eu preciso de alguma ajuda, ou sei lá o quê... Mas nunca falar uma coisa dessas!
			
			-- Eu falei isso pra ela. Falei que, ao invés dela falar isso, talvez era ela que devia ver o que `tava' acontecendo, já que você nunca teve problemas antes e só `tá' com problemas na matéria dela. Ela ficou pálida, pediu desculpas dizendo que ela não pode fazer isso porque `tá' sem tempo, disse que vai viajar pra outra cidade, que vai mudar de escola... A única coisa que ela me perguntou foi se eu não ficava preocupada com você faltando. Eu dei risada. \emph{\textquotedblleft Eu nunca, na minha vida, precisei obrigar a Fernanda a `vim' pra escola. Ela vem porque ela quer. Quando ela não quer `vim', ela não vem. E, se ela está dizendo que quer ficar em casa pra resolver essa lista, eu confio nela, ela sabe o que `tá' fazendo. Além do mais, ela pode até não `tá' vindo pra escola, mas `tá' lá em casa, estudando. Ela é inteligente, o resto das matérias ela pega depois.\textquotedblright }
			
			-- Eu vou resolver essa lista! Ainda mais depois disso que ela disse! É uma questão de honra agora!
			
			-- Eu sei que você vai, você sempre consegue!
		\end{quote}
		E foi assim que eu passei a semana toda, dentro de casa, diante daquela lista e de um livro de Matemática que parecia não me ajudar em nada. Até que, numa tarde que estava bastante tediosa, eu notei que, sobre a mesa onde eu estudava, havia um espelho de mão. E foi justamente no meio de uma brincadeira, que eu resolvi fazer com aquele espelho para desanuviar um pouco, que eu acabei percebendo uma coisa interessante: ao enxergar a imagem que o espelho construiu de uma régua que também estava sobre aquela mesa, eu percebi que aquela imagem era \textbf{revertida} conforme, por exemplo, a Figura \ref{regua-e-imagem} está ilustrando.
		E como colocar o espelho perpendicularmente àquela régua, no seu ponto $ \boldsymbol{0} $, me fez perceber que a justaposição dessa régua com a sua imagem correspondia a um trecho da reta dos números reais, desde que fosse feita a identificação de cada número \reflectbox{$ \boldsymbol{Z} $} (ou seja, cada número $ \boldsymbol{Z} $ revertido que aparece na imagem) como $ - \boldsymbol{Z} $, foi exatamente isso que permitiu com que eu justificasse todas aquelas famigeradas regras de sinais que estavam envolvidas para com os números positivos e negativos. Como?
		\begin{figure}[!t]
			\centering
			\includegraphics[viewport=360 10 0 230,scale=1.3]{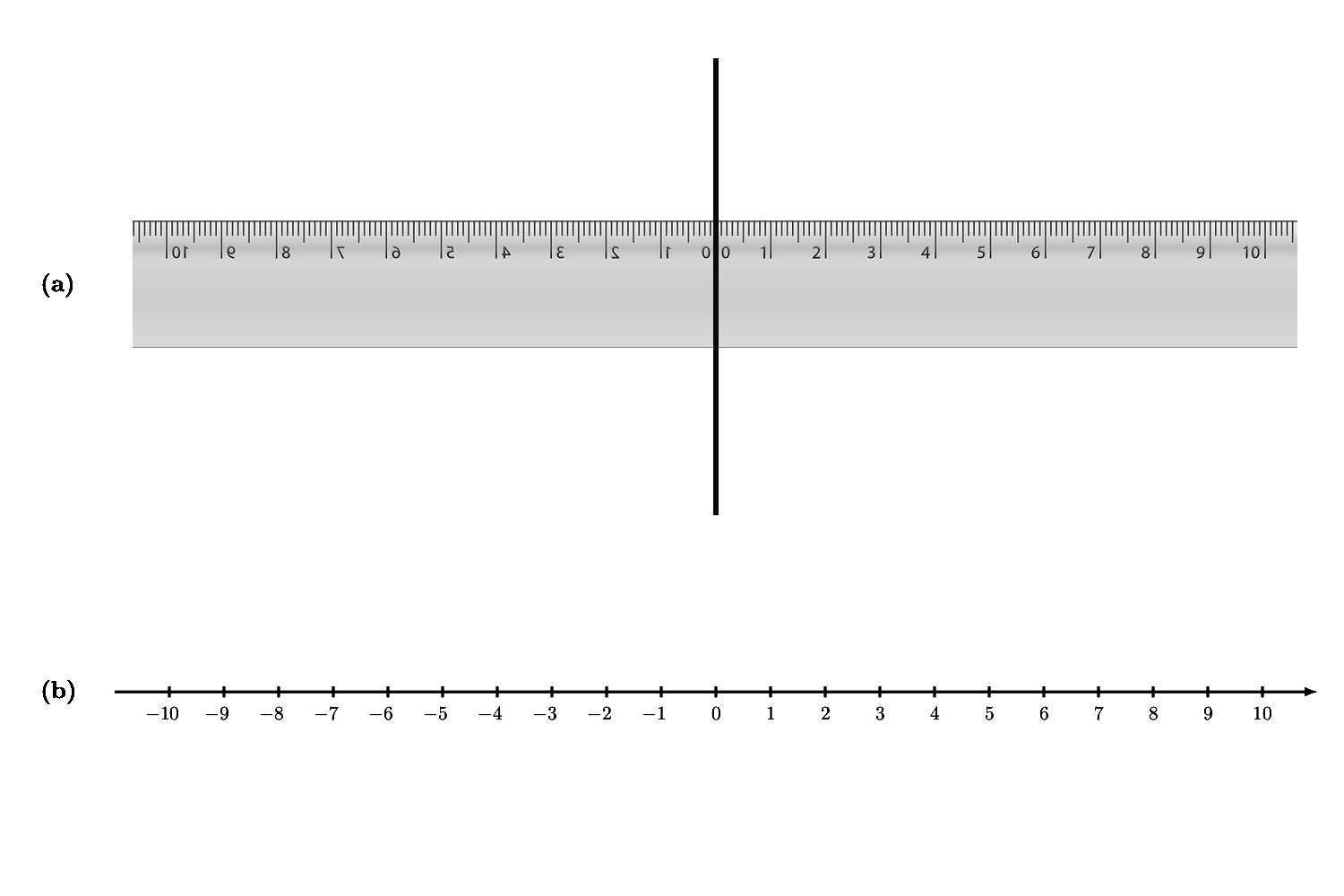}
			\caption{\label{regua-e-reta}\textbf{(a)} Justaposição de uma régua (à direita) com a sua imagem (à esquerda) mediante a alocação de um espelho plano (que está destacado apenas em preto, sem dar qualquer destaque à sua superfície refletora apenas para realçar a simetria desta figura) perpendicularmente a essa régua no seu ponto $ \boldsymbol{0} $. \textbf{(b)} Reta dos números reais onde está sendo dado destaque aos vinte números inteiros não nulos que estão mais próximos de $ \boldsymbol{0} $.} 
		\end{figure}
			
		A primeira coisa que eu percebi (diante dessa justaposição que está sendo bem ilustrada pela Figura \ref{regua-e-reta}(a)) foi que, do mesmo jeito que o número $ \boldsymbol{5} $ estava mais distante do número $ \boldsymbol{0} $ que o número $ \boldsymbol{3} $, o número \reflectbox{$ \boldsymbol{5} $} também estava mais distante de $ \boldsymbol{0} $ que o número \reflectbox{$ \boldsymbol{3} $}. E isso só acontecia porque todas as distâncias que existiam entre os números da régua real eram \textbf{preservadas} pelos da régua revertida.
		\begin{figure}[!t]
			\centering
			\includegraphics[viewport=360 10 0 305,scale=1.3]{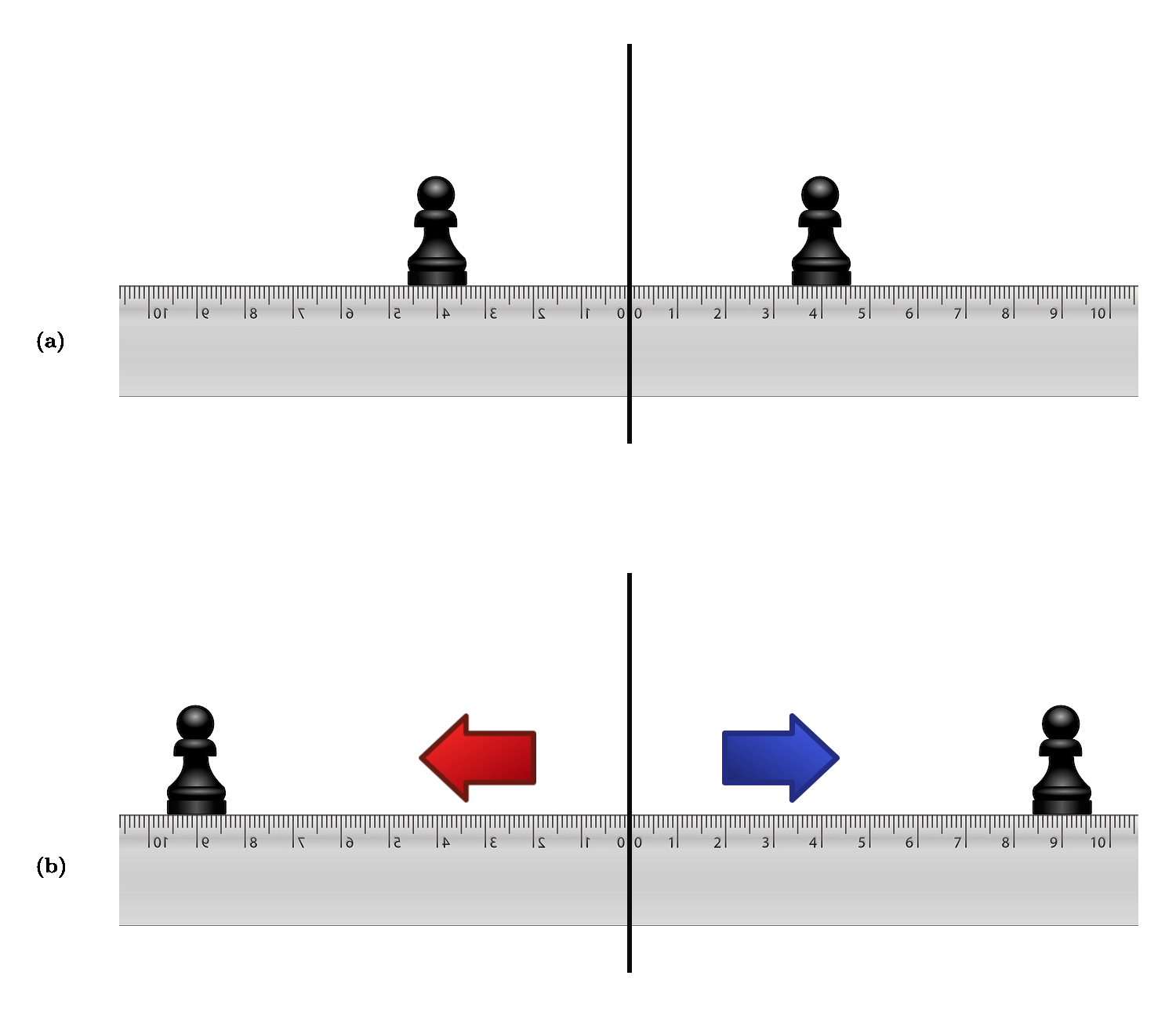}
			\caption{\label{soma-total}\textbf{(a)} Situação inicial onde um objeto (o peão preto) foi colocado sobre o número $ \boldsymbol{4} $ da régua (à direita) e a sua imagem (que foi construída pelo espelho plano que está novamente destacado em preto, sem qualquer destaque à sua superfície refletora) está sobre o número \quatrorevertido . \textbf{(b)} Situação final onde o mesmo objeto e, consequentemente, a sua imagem foram deslocados para os números $ \boldsymbol{9} $ e \noverevertido \hspace*{0.01cm} respectivamente, caracterizando ambos os deslocamentos como uma operação de soma já que tanto o objeto como a sua imagem estão mais distantes de $ \boldsymbol{0} $.}
		\end{figure}
			
		Já a segunda coisa que eu percebi foi que, quando eu colocava um objeto (como, por exemplo, o peão que consta na Figura \ref{soma-total}) sobre a posição $ \boldsymbol{4} $ da régua e o deslocava até a posição $ \boldsymbol{9} $, a imagem daquele objeto também saía da posição \reflectbox{$ \boldsymbol{4} $} e chegava até a posição \reflectbox{$ \boldsymbol{9} $}. E como esse processo de sair da posição $ \boldsymbol{4} $ da régua e chegar até a posição $ \boldsymbol{9} $ pode ser interpretado em termos da \textbf{soma}
		\begin{equation*}
			\boldsymbol{4} + \boldsymbol{5} = \boldsymbol{9} \ ,
		\end{equation*}
		não havia como não interpretar o processo de sair da posição \reflectbox{$ \boldsymbol{4} $} e chegar até a posição \reflectbox{$ \boldsymbol{9} $} também em termos de uma soma
		\begin{equation*}
			\quatrorevertido + \cincorevertido = \noverevertido \ ,
		\end{equation*}
		já que ambos os processos levavam a números que mantinham a mesma distância (de nove unidades) em relação ao número $ \boldsymbol{0} $.
		\begin{figure}[!t]
			\centering
			\includegraphics[viewport=360 10 0 305,scale=1.3]{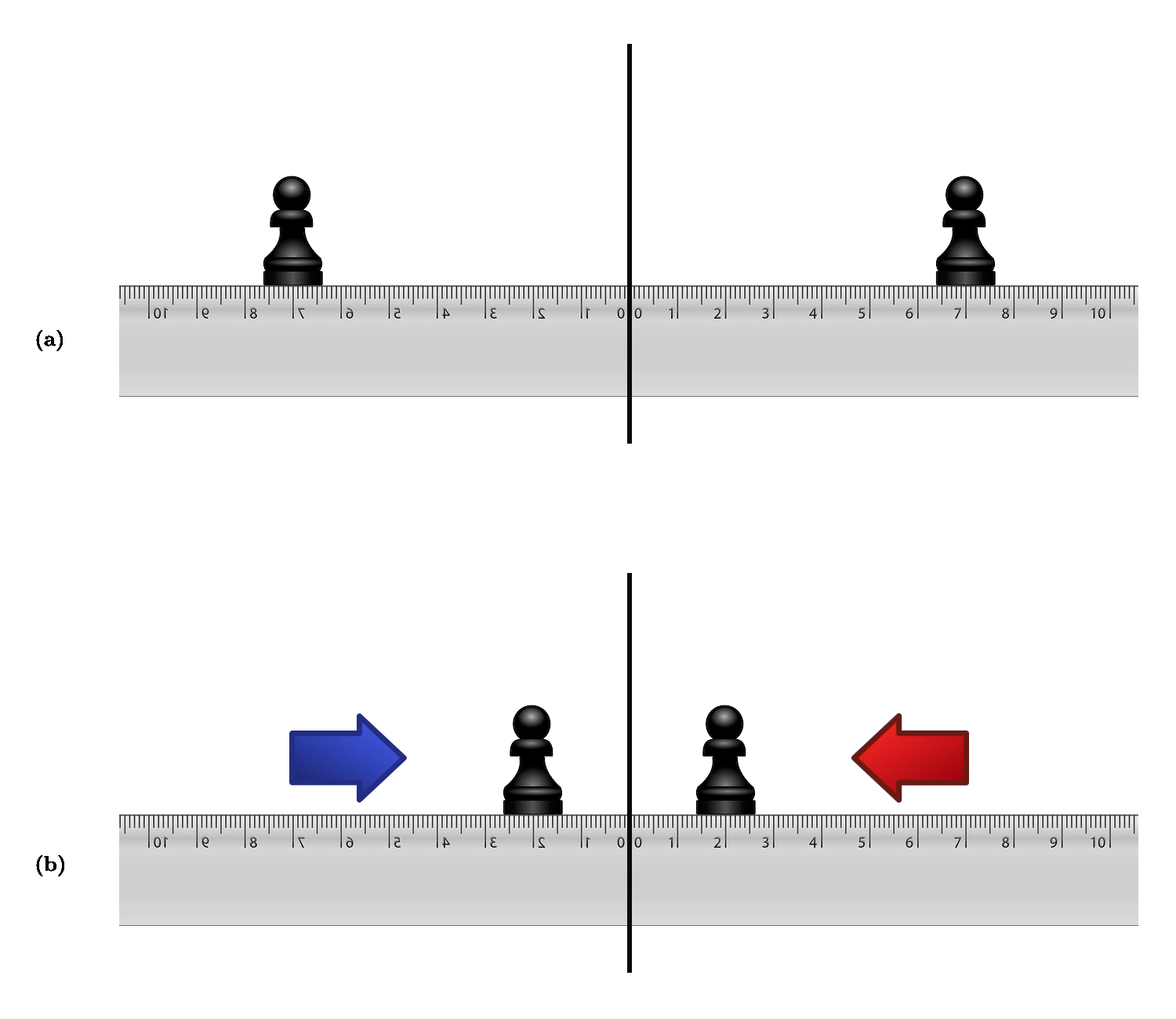}
			\caption{\label{subtracao-total}\textbf{(a)} Nova situação inicial onde o mesmo objeto da figura anterior (que continua à direita) foi colocado sobre o número $ \boldsymbol{7} $ da régua e, consequentemente, a sua imagem (que foi construída pelo mesmo espelho plano da figura anterior) agora está sobre o número \seterevertido . \textbf{(b)} Nova situação final onde tal objeto e, consequentemente, a sua imagem foram deslocados para os números $ \boldsymbol{2} $ e \doisrevertido \hspace*{0.01cm} respectivamente, o que caracteriza ambos os deslocamentos como uma operação de subtração já que tanto o objeto como a sua imagem agora estão mais próximos de $ \boldsymbol{0} $.}
		\end{figure}
			
		Seguindo essa mesma linha de raciocínio, a terceira coisa que eu observei foi que, como o processo de sair da posição $ \boldsymbol{7} $ da régua e chegar até a posição $ \boldsymbol{2} $ pode ser interpretado em termos da \textbf{subtração}
		\begin{equation*}
			\boldsymbol{7} - \boldsymbol{5} = \boldsymbol{2} \ ,
		\end{equation*}
		também não havia como eu não interpretar o processo de sair da posição \reflectbox{$ \boldsymbol{7} $} e chegar até a posição \reflectbox{$ \boldsymbol{2} $} em termos de uma subtração
		\begin{equation*}
			\seterevertido - \cincorevertido = \doisrevertido \ ,
		\end{equation*}
		já que ambos os processos levavam a números que mantinham a mesma distância (de duas unidades) em relação ao número ponto $ \boldsymbol{0} $.
			
		Diante de tudo isso, como, pra mim, estava bastante clara a existência de uma correspondência entre a reta dos números reais e uma outra reta, que poderia ser fisicamente construída com o auxílio de uma régua semi-infinita desde que
		\begin{itemize}
			\item[\textbf{(i)}] um espelho plano fosse posicionado \textbf{perpendicularmente} a essa régua sobre o seu ponto $ \boldsymbol{0} $ e
			\item[\textbf{(ii)}] todo número \reflectbox{$ \boldsymbol{Z} $} fosse identificado como $ \boldsymbol{-Z} $,
		\end{itemize}
		todas essas regras de sinais ficaram justificadas. Afinal, como o processo que pode ser interpretado como uma soma (o da Figura \ref{soma-total}) está mostrando que a soma de um número negativo está associada a um deslocamento para a esquerda (o qual está sendo destacado pela flecha vermelha da Figura \ref{soma-total}\textbf{(b)}), tal como também já ocorre com a subtração de um número positivo (o que também está propositalmente destacado por uma flecha vermelha na Figura \ref{subtracao-total}\textbf{(b)}), é possível afirmar que
		\begin{equation*}
			- \left( + \boldsymbol{Z} \right) = + \left( - \boldsymbol{Z} \right) = - \boldsymbol{Z} \ .
		\end{equation*}
		Analogamente, como o deslocamento que caracteriza a subtração de um número negativo é para a direita (conforme bem ilustra a Figura \ref{subtracao-total}\textbf{(b)} com o auxílio de uma flecha azul), assim como também acontece com a soma de um número positivo (o que também está sendo propositalmente destacado, pela flecha de mesma cor, na Figura \ref{soma-total}\textbf{(b)}), também passa a ser possível afirmar que
		\begin{equation*}
			- \left( - \boldsymbol{Z} \right) = + \left( + \boldsymbol{Z} \right) = + \boldsymbol{Z} \ .
		\end{equation*}
		
		Depois que eu enxerguei toda essa correspondência, tudo se tornou tão trivial que eu consegui fazer todas as contas daquela lista numa única tarde. E, como eu morava bem perto da escola onde eu estudava, assim que terminei eu já fui até lá, no início da noite, entregar a lista resolvida para a minha professora antes que as suas aulas terminassem. Cerca de duas ou três semanas depois, com todas as listas já corrigidas, a professora começa uma chamada anunciando, em voz alta, as notas de todos.
		\begin{quote}
			-- Maria Fernanda.
			
			-- Presente.
			
			-- \textbf{Dez}, meus parabéns!
		\end{quote}
		Naquele momento, todos na sala de aula se viraram para mim num ato único, olhando atônitos. Deu até para ouvir o barulho de todos se virando ao mesmo tempo. Afinal de contas, depois de eu ter ficado de recuperação em Matemática e de ter sofrido \textquotedblleft bullying\textquotedblright \hspace*{0.01cm} ao longo de todo o ano justamente daquele grupo onde o sorteio me inseriu, eu larguei esse grupo (coisa que ninguém, ali, sabia até então já que eu desapareci de todas as aulas), resolvi toda aquela lista de exercícios sozinha, sem a ajuda de ninguém, e ainda fui a única pessoa que tirou \textbf{dez} não apenas naquela lista, mas naquele bimestre.
		\begin{quote}
			-- O que foi, nunca me viram?
		\end{quote}
		Depois disso, eu voltei ao meu normal e nunca mais tirei uma nota baixa em Matemática no Ensino Básico, do mesmo jeito que eu nunca mais vi ninguém daquela sala praticar \textquotedblleft bullying\textquotedblright \hspace*{0.01cm} contra mim nem contra nenhum dos meus amigos ao meu redor.
			
	\section{Mas por que usar um espelho plano?}
	
		Na verdade, o fato de eu ter usado um espelho plano nesse experimento, que me fez entender todas essas regras de sinais, não foi uma questão de escolha: aquele era o único espelho que eu tinha ao meu alcance. Entretanto, é importante deixar bem claro que são justamente os espelhos planos que são os mais indicados para este fim. E para deixar bem claro por que eles são realmente os mais indicados, é interessante nós voltarmos as nossas atenções para \textquotedblleft outros\textquotedblright \hspace*{0.01cm} tipos de espelhos cujo nome parece não remeter à nada que é plano: os \textbf{espelhos esféricos}.
		
		\subsection{\label{subsection:paraxiais}O que são espelhos esféricos?}
		
\begin{figure}[!t]
				\centering
				\includegraphics[viewport=460 10 0 305,scale=1.0]{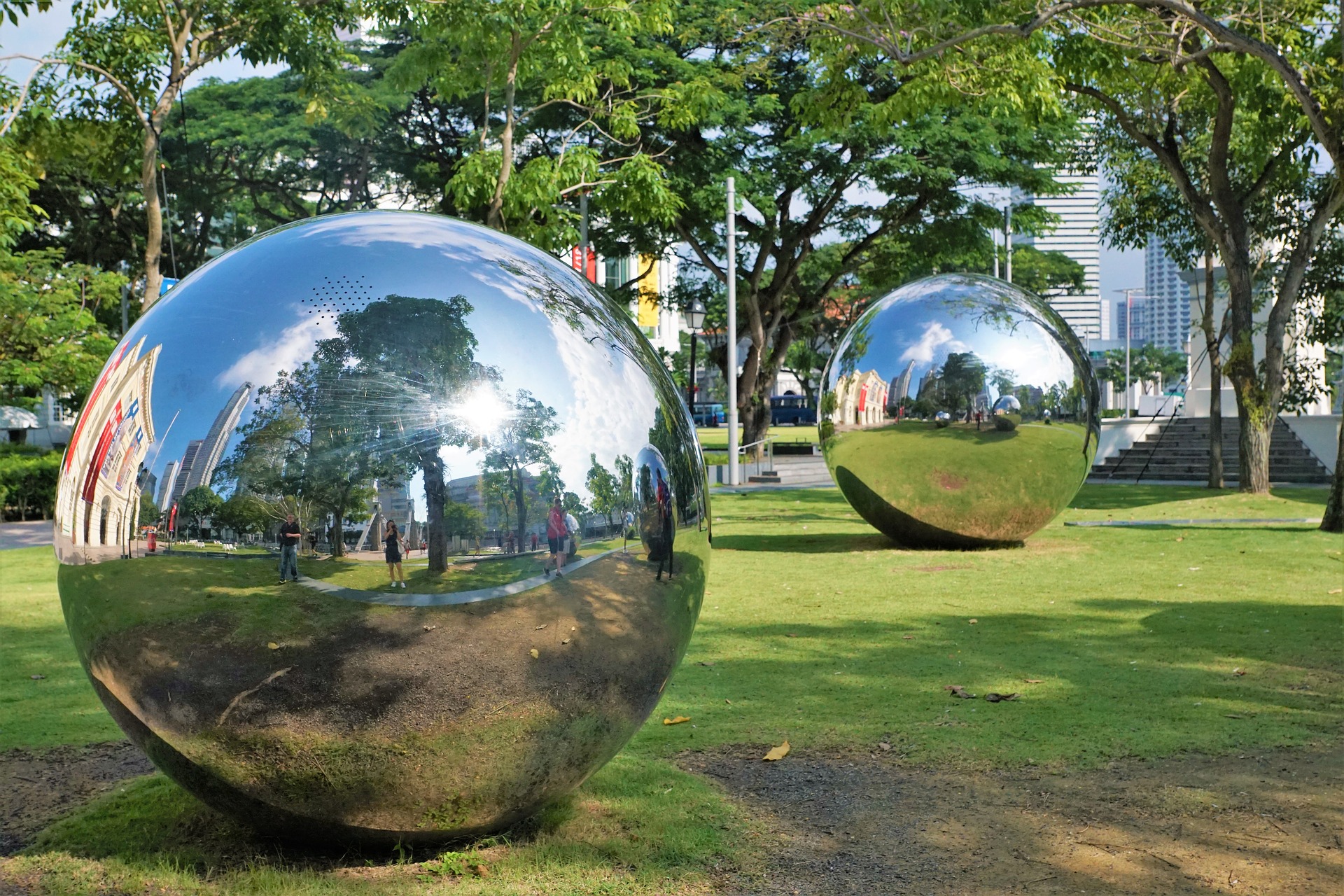}
				\caption{\label{mirror-balls}Duas esferas cujas superfícies exteriores foram preparadas para funcionarem como espelhos. Note que, para quem olha essas duas esferas frontalmente, as imagens que aparecem sobre as suas superfícies até que são bastante nítidas na região mais central, embora sejam um pouco deformadas nas outras regiões devido à não planitude dessas superfícies.}
			\end{figure}

			De um modo geral, é possível afirmar que os espelhos esféricos recebem esse nome por serem fabricados com o auxílio de alguma estrutura esférica. Aliás, a Figura \ref{mirror-balls} já traz um bom exemplo de tais espelhos, uma vez que as duas esferas que nela consta possuem superfícies que conseguem formar imagens até que bem nítidas de diversos objetos que estão nas suas adjacências.
			
\begin{figure}[!t]
				\centering
				\includegraphics[viewport=460 10 0 460,scale=1.0]{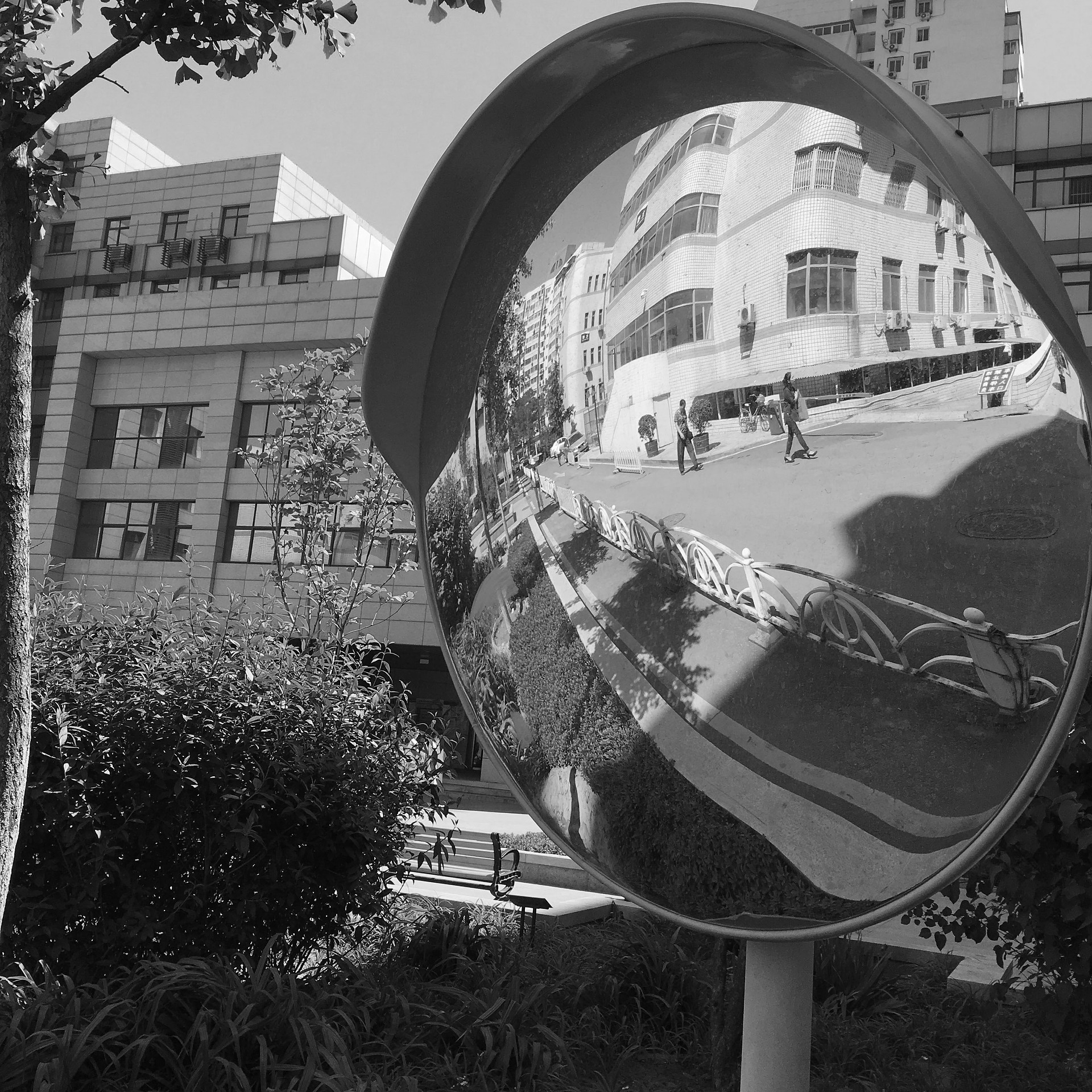}
				\caption{\label{spherical-mirror}Espelho esférico comumente usado nas situações onde se faz necessário diminuir os tamanhos das imagens dos objetos que estão próximos a ele aumentando, assim, o seu campo visual.}
			\end{figure}

\begin{figure}[!t]
				\centering
				\includegraphics[viewport=180 10 0 170,scale=1.1]{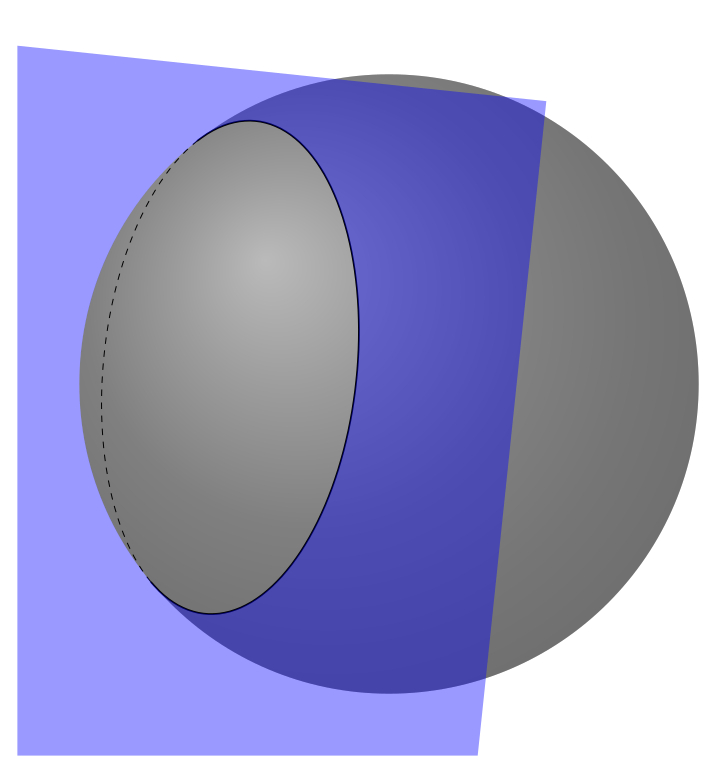}
				\caption{\label{corte-esfera}Esquema que ilustra a lógica por trás do corte que é feito numa esfera oca (destacada em cinza) com o auxílio de um plano (destacado com um azul semitransparente) para formar uma calota esférica que pode dar suporte a um espelho esférico como, por exemplo, o que aparece na Figura \ref{spherical-mirror}.}
			\end{figure}

			É claro que vale notar que nem todo espelho esférico possui essa forma explicitamente esférica dos que aparecem nessa Figura \ref{mirror-balls} e um bom exemplo disso é o espelho esférico que aparece na Figura \ref{spherical-mirror}.
			No entanto, apesar desse último espelho não ter um formato explicitamente esférico, o que justifica o seu nome é o fato dele possuir o formato de uma \textbf{calota esférica}: ou seja, conforme bem ilustra a Figura \ref{corte-esfera}, o seu processo de fabricação pode ser visto em termos de um corte numa esfera oca, a qual teve (pelo menos) uma das suas superfícies preparadas para refletir luz e formar imagens nítidas.
		
\begin{figure}[!t]
				\centering
				\includegraphics[viewport=330 10 0 180,scale=1.3]{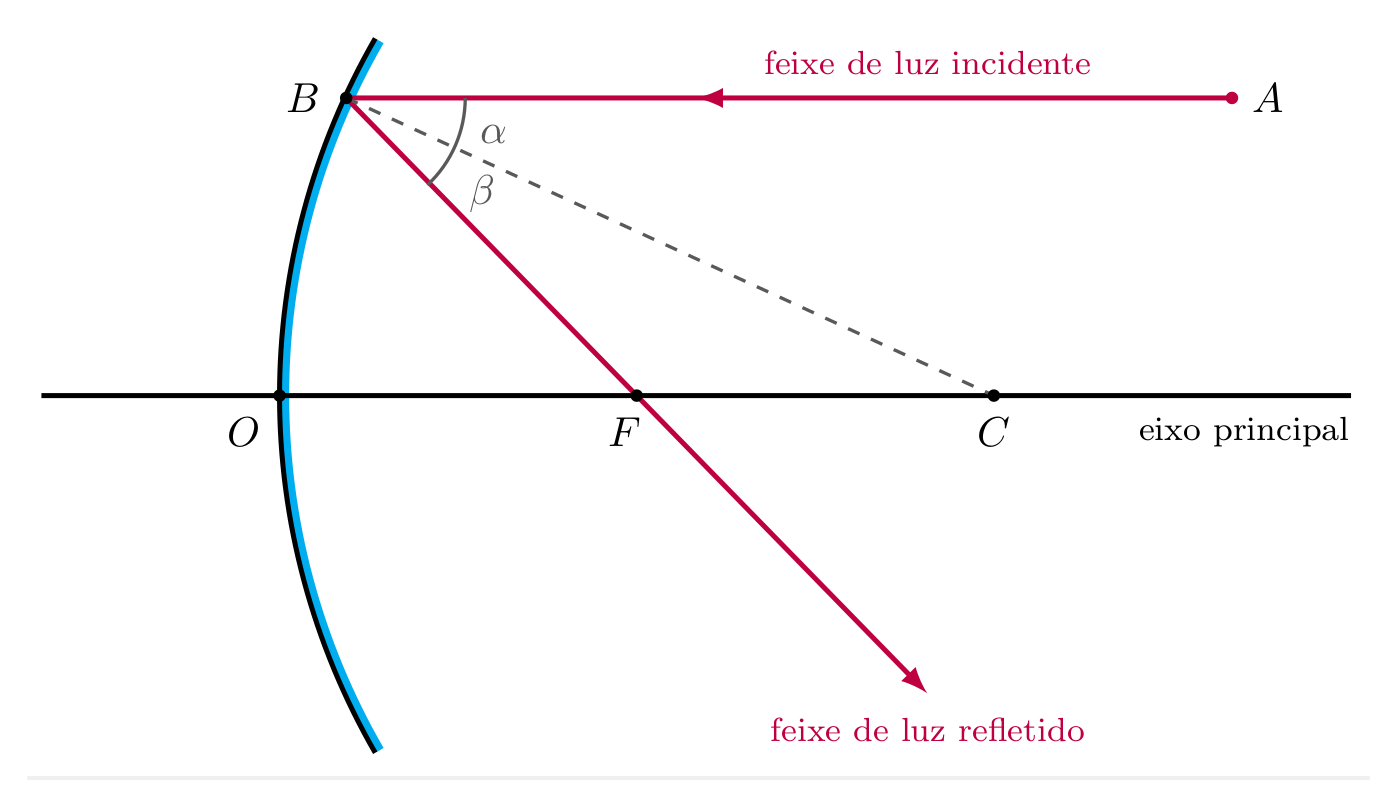}
				\caption{\label{light-reflection}Corte lateral de um espelho \textbf{côncavo} (ou seja, de um espelho cuja face espelhada (destacada em ciano) é a superfície interna da calota que lhe deu origem) onde é dado destaque ao seu eixo principal que o intercepta no centro $ \boldsymbol{O} $ da sua superfície. Note que o feixe de luz incidente está sendo refletido para o ponto $ \boldsymbol{F} $ (o foco) deste espelho esférico, ponto o qual não coincide com o centro $ \boldsymbol{C} $ da esfera oca que deu origem a este espelho. Como tal reflexão sempre ocorre com $ \boldsymbol{\alpha } = \boldsymbol{\beta } $, essa focalização de todos os feixes em $ \boldsymbol{F} $ só ocorre quando o ângulo $ \boldsymbol{\measuredangle BCO} $ é extremamente pequeno (ou seja, quando todos esses feixes são paraxiais \cite{sutanto}).}
			\end{figure}

			Diante desta última observação, e independentemente de um espelho ter sido construído usando uma esfera completa ou apenas um pedaço dela, o fato é que todo espelho esférico pode ser caracterizado por um raio $ \boldsymbol{R} $ -- no caso, pelo raio da esfera que lhe deu origem. E uma relação, que é bem conhecida entre $ \boldsymbol{R} $ e as posições do objeto $ \boldsymbol{p_{\mathrm{ob}}} $ e da imagem $ \boldsymbol{p_{\mathrm{im}}} $ que esse espelho forma do objeto, é
			\begin{equation}
				\frac{\boldsymbol{1}}{\boldsymbol{p_{\mathrm{ob}}}} + \frac{\boldsymbol{1}}{\boldsymbol{p_{\mathrm{im}}}} = \frac{\boldsymbol{2}}{\boldsymbol{R}} \ . \label{gauss-1}
			\end{equation}
			Esta é a chamada \textbf{equação de Gauss}, a qual costuma ser popularmente expressa como
			\begin{equation}
				\frac{\boldsymbol{1}}{\boldsymbol{p_{\mathrm{ob}}}} + \frac{\boldsymbol{1}}{\boldsymbol{p_{\mathrm{im}}}} = \frac{\boldsymbol{1}}{\boldsymbol{f}} \ , \label{gauss-2}
			\end{equation}
			onde a razão
			\begin{equation}
				\boldsymbol{f} = \boldsymbol{R} / \boldsymbol{2} \label{d-focal}
			\end{equation}
			é reconhecida como a \textbf{distância focal} desses espelhos esféricos. Essa razão só é assim chamada por ser a distância que existe entre a superfície de um espelho esférico e o ponto para o qual todos os feixes de luz \textbf{paraxiais} (ou seja, que incidem sobre essa superfície paralelamente ao eixo de simetria e muito próximos a ele) convergem. Note que essa é justamente a situação que está sendo ilustrada na Figura \ref{light-reflection}, onde nós vemos um único feixe de luz partindo do ponto $ \boldsymbol{A} $ e incidindo sobre a superfície interna de uma calota esférica: se essa superfície tiver sido bem preparada para funcionar como um espelho, esse feixe será refletido de um ponto $ \boldsymbol{B} $ dessa calota para um ponto $ \boldsymbol{F} $ que pertence ao eixo de simetria deste espelho (eixo esse que será chamado, a partir de agora, como \textbf{eixo principal}). E como todos os feixes de luz paraxiais sempre serão refletidos para esse mesmo ponto $ \boldsymbol{F} $ é justamente esse ponto que pode ser reconhecido como o \textbf{foco} deste espelho, o qual se mantém a uma distância $ \boldsymbol{f} $ da superfície em questão.
		
		\subsection{Qual é a relação que existe entre espelhos planos e esféricos?}
		
			Embora a demostração da validade de todas essas relações (\ref{gauss-1}), (\ref{gauss-2}) e (\ref{d-focal}) possa ser encontrada em diversas fontes bibliográficas \cite{hecht,savelyev,moyses}, uma demonstração bastante didática também consta no Apêndice \ref{appendix} deste artigo por uma questão de completude. Todavia, diante da proposta deste artigo, a pergunta que não quer calar é: por que será que é tão importante nós termos consciência de todas essas informações sobre os espelhos esféricos? Aliás, será que existe alguma relação entre esses espelhos e os espelhos planos, a qual consegue explicar por que os últimos são os mais indicados nesta tarefa de entender as regras de sinais da Matemática?
			
			Para entender as respostas de todas essas questões é importante notar, primeiro, que, às vezes, quando nós estamos diante de uma esfera que possui um raio gigantesco (quando comparado, por exemplo, às nossas dimensões), fica muito difícil reconhecer que nós estamos realmente diante de uma esfera. Esse é o caso, por exemplo, da Terra: dependendo da região onde nós estamos sobre ela, é impossível reconhecer o seu formato esférico apenas visualmente. Agora imagine uma situação onde o raio de uma esfera tem um tamanho que tende ao infinito: neste caso é \textbf{impossível} distinguir a superfície de uma esfera da de um plano infinito. E como, neste caso onde o raio de uma esfera tem um tamanho que tende ao infinito, a expressão (\ref{d-focal}) nos mostra que
			\begin{equation}
				\lim _{\boldsymbol{R} \rightarrow \boldsymbol{\infty }} \left( \frac{\boldsymbol{1}}{\boldsymbol{f}} \right) = \boldsymbol{0} \ , \label{limite}
			\end{equation}
			a substituição deste resultado em (\ref{gauss-2}) nos leva a
			\begin{equation}
				\boldsymbol{p_{\mathrm{im}}} = - \boldsymbol{p_{\mathrm{ob}}} \ . \label{resultado}
			\end{equation}
			
\begin{figure}[!t]
				\centering
				\includegraphics[viewport=340 10 0 170,scale=1.4]{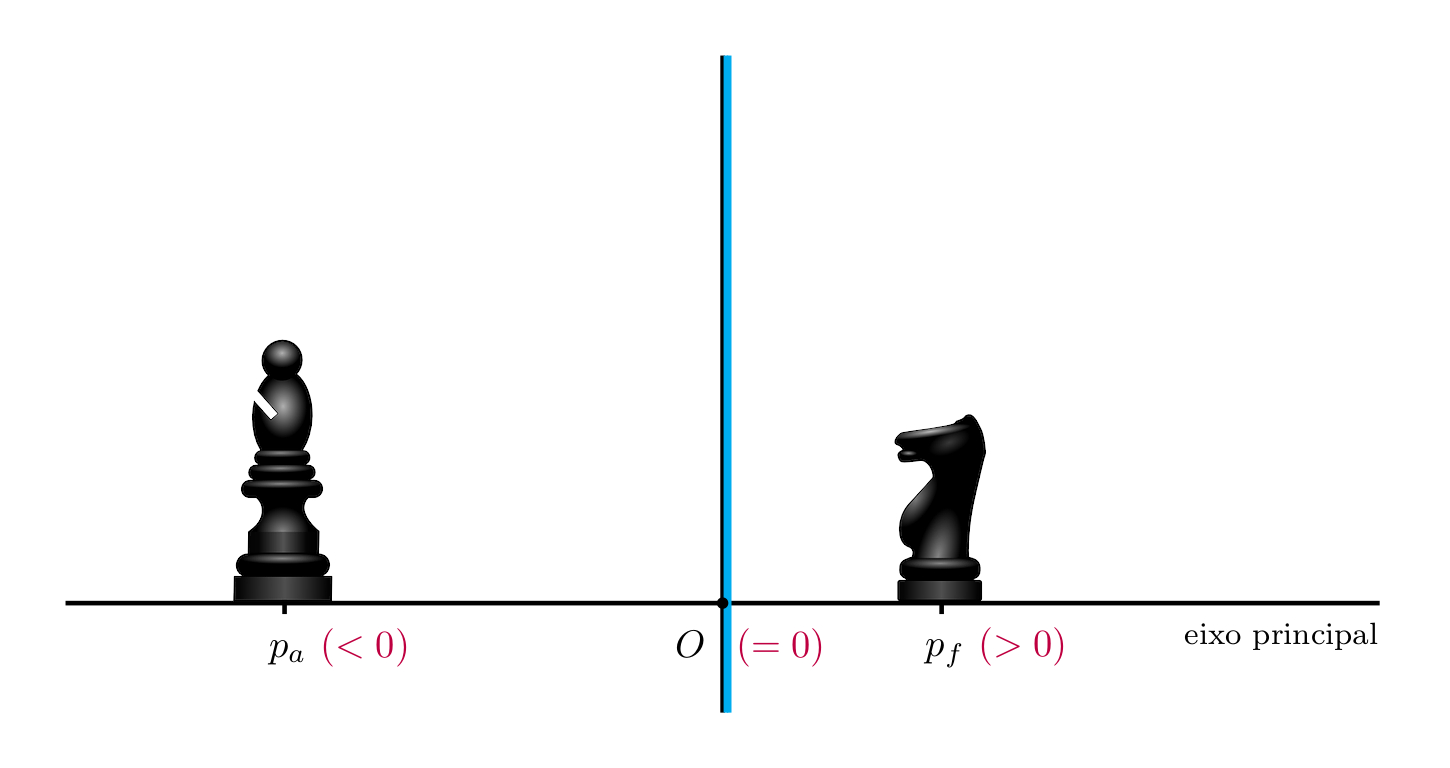}
				\caption{\label{convencao}Convenção que é adotada para descrever todas as posições que estão relacionadas à formação de imagens de objetos por um espelho esférico: sempre que um objeto e/ou a sua imagem estiverem diante ou atrás da face espelhada deste espelho, as suas posições serão representadas por números reais positivos ou negativos respectivamente. Este é justamente o caso do objeto (o cavalo preto) que aparece à direita: como a face espelhada está voltada para esse lado direito, a sua posição $ \boldsymbol{p_{f}} $ é descrita por um número real positivo, o qual pode ser perfeitamente interpretado como a distância que existe entre tal objeto e o espelho. Já no caso do outro objeto (o bispo preto) que consta à esquerda (mais especificamente atrás da face espelhada), a sua posição $ \boldsymbol{p_{a}} $ é descrita por um número real negativo. Este número negativo também pode ser interpretado em termos da distância que existe entre esse último objeto e o espelho, só que multiplicada por $ - \boldsymbol{1} $. Aqui, vale notar que os objetos que constam dos dois lados desta figura foram escolhidos apenas a título de ilustração e, portanto, não devem ser interpretados como a imagem um do outro. Outra coisa que também vale notar aqui é que, embora o espelho que aparece nesta figura seja plano, tudo o que foi dito aqui também vale para um espelho esférico, até mesmo porque espelhos planos são interpretáveis como espelhos esféricos que possuem raios infinitos.}
			\end{figure}

			Já a segunda coisa que nós podemos fazer para entender as respostas dessas questões é notar que, quando nós estamos diante da equação de Gauss (esteja ela expressa sob a forma (\ref{gauss-1}) ou sob a forma (\ref{gauss-2})), as posições $ \boldsymbol{p_{\mathrm{ob}}} $ de um objeto e $ \boldsymbol{p_{\mathrm{im}}} $ da imagem que um espelho esférico forma para esse objeto acabam sendo localizadas com o auxílio do eixo principal. Ou seja, conforme a Figura \ref{convencao} está ilustrando,
			se nós colocarmos um objeto diante de um espelho esférico, tais posições são localizadas com o auxílio da métrica que é atribuída a esse eixo, a qual, por uma questão de \textquotedblleft convenção\textquotedblright , está em correspondência para com a reta dos números reais. Afinal, de acordo com o que essa Figura \ref{convencao} nos mostra, enquanto todas as posições na frente do espelho são descritas por números positivos, todas as que ficam atrás do espelho são descritas por números negativos. Aliás, devido a essa correspondência, vale notar que a origem $ \boldsymbol{0} $ desta métrica está sobre o ponto $ \boldsymbol{O} $ onde o espelho e o seu eixo principal se interceptam.
			
\begin{figure}[!t]
				\centering
				\includegraphics[viewport=330 10 0 155,scale=1.5]{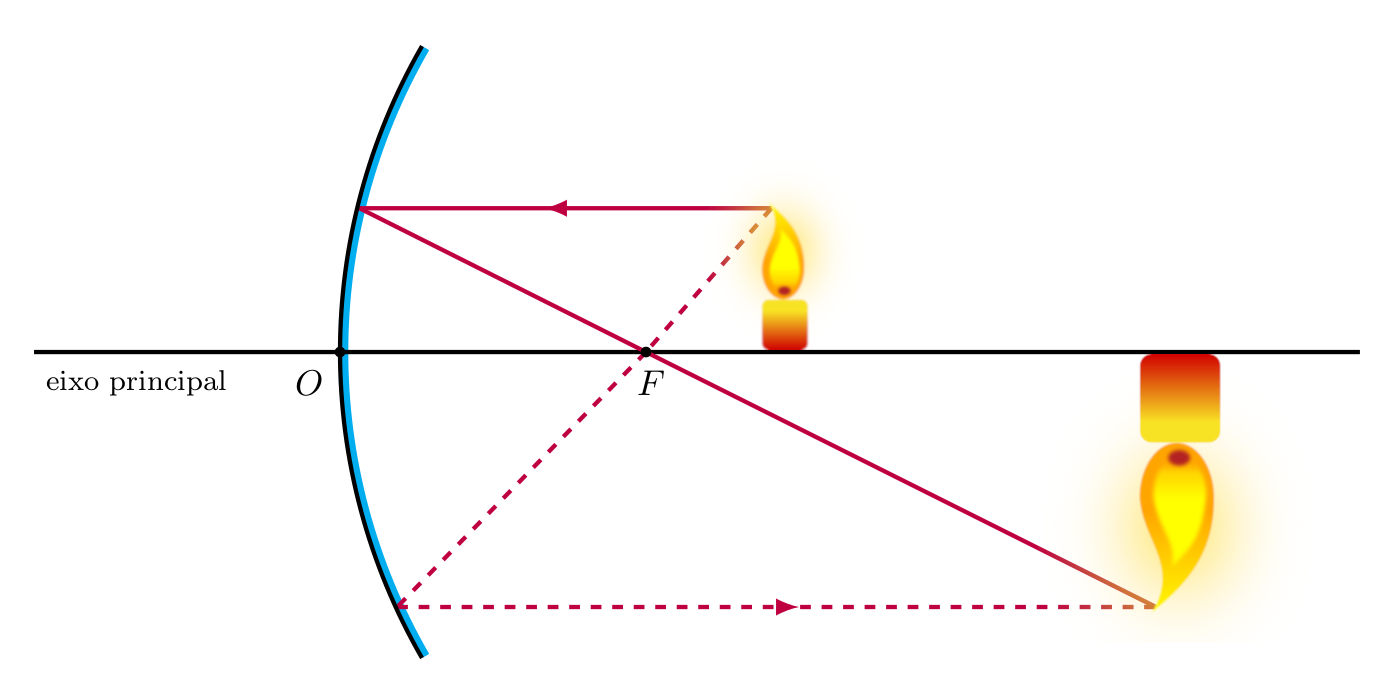}
				\caption{\label{imagem-real}Corte lateral de um espelho \textbf{côncavo} que consegue ilustrar como uma imagem real pode ser formada. Para isso, um objeto real (uma vela) foi colocado na frente desse espelho, bem próximo ao seu foco, e a sua imagem real (por também estar na frente desse espelho) pôde ser identificada numa posição um pouco mais distante. As duas linhas (uma contínua e outra hachurada) que estão destacadas nesta figura (com a cor púrpura) representam dois feixes de luz que partem do objeto real em direções distintas, os quais voltam a se cruzar no lugar onde a imagem é identificada. Maiores detalhes sobre o porquê tais feixes foram destacados podem ser encontrados, por exemplo, no Apêndice \ref{appendix}.}
			\end{figure}

			Em linhas bem gerais, é válido afirmar que um dos motivos que justificam essa \textquotedblleft convenção\textquotedblright \hspace*{0.01cm} é, por exemplo, o fato de espelhos esféricos serem capazes de formar dois tipos de imagens. E uma delas é a que se pode chamar de \textbf{imagem virtual}, a qual só recebe este nome por só poder ser enxergada dentro do espelho: esse é justamente o caso, por exemplo, de todas as imagens que já puderam ser vistas nas Figuras \ref{mirror-balls} e \ref{spherical-mirror}. Já o outro tipo de imagem que um espelho esférico é capaz de formar é a \textbf{imagem real}, a qual recebe este nome por poder ser vista fora do espelho (isto é, no \textquotedblleft mundo real\textquotedblright , bem ao nosso lado) conforme a Figura \ref{imagem-real} está ilustrando.
			Ou seja, por mais que seja intuitivo interpretar $ \boldsymbol{p_{\mathrm{ob}}} $ como a distância que existe entre um objeto real e um espelho (já que é justamente essa distância que nós conseguimos medir na prática), o fato da imagem de um objeto poder (como uma possibilidade) ser formada tanto do lado de dentro como do lado de fora de um espelho parece nos impelir a dar uma interpretação um pouco diferente à posição $ \boldsymbol{p_{\mathrm{im}}} $, já que $ \boldsymbol{p_{\mathrm{im}}} $ precisa ser capaz de fazer tal distinção. Nestes termos, como atribuir tal escala métrica ao eixo principal nos permite não apenas distinguir as posições de imagens reais e virtuais usando números reais positivos e negativos, mas também a interpretar $ \left\vert \boldsymbol{p_{\mathrm{im}}} \right\vert $ (ou seja, o módulo ou a norma de $ \boldsymbol{p_{\mathrm{im}}} $) como a distância que essa imagem (real ou não) mantém do espelho, o uso dessa escala acaba sendo muito bem-vindo.
			
			Diante de todas as observações que acabaram de ser feitas aqui, e principalmente diante do resultado que foi obtido em (\ref{resultado}), não é nada difícil concluir que, se nós construírmos qualquer reta perpendicular a um espelho plano num ponto $ \boldsymbol{O} $ arbitrário e darmos, à ela, uma métrica que tenha a sua origem em $ \boldsymbol{O} $, essa reta se estenderá para dentro do espelho. E como essa métrica se estenderá para dentro desse espelho sem sofrer qualquer tipo de deformação devido ao resultado (\ref{limite}), também não será nada difícil constatar que, para cada ponto $ \boldsymbol{p_{\mathrm{ob}}} $ que existir sobre a reta real (ou seja, sobre a reta que nós realmente construímos diante da superfície espelhada), a sua imagem estará, atrás do espelho, sobre a posição $ \boldsymbol{p_{\mathrm{im}}} = - \boldsymbol{p_{\mathrm{ob}}} $ que é identificável sobre essa reta estendida. 
	
	\section{Atividades que podem ser propostas}
	
		Se nós analisarmos esta situação desta construção pelo ponto de vista lúdico, o ato de colocar uma régua, ou qualquer outro objeto que possa ser interpretado como uma reta numérica, bem na frente da superfície espelhada de um espelho plano, abre uma espécie de \textquotedblleft portal\textquotedblright \hspace*{0.01cm} através do qual é possível enxergar o \textquotedblleft habitat\textquotedblright \hspace*{0.01cm} dos números negativos. Afinal, note que, conforme já bem ilustrou a Figura \ref{regua-e-reta}\textbf{(a)}, para qualquer número $ \boldsymbol{Z} $ que venha a ser escolhido sobre a régua que está na frente do espelho, sempre será possível enxergar o seu oposto \reflectbox{$ \boldsymbol{Z} $}, dentro do espelho, à mesma distância da sua superfície. Assim, como o reconhecimento de \reflectbox{$ \boldsymbol{Z} $} como $ - \boldsymbol{Z} $ entra em pleno acordo com a propriedade
		\begin{equation}
			\boldsymbol{Z} + \left( - \boldsymbol{Z} \right) = - \boldsymbol{Z} +  \boldsymbol{Z} = \boldsymbol{0} \label{reais-1}
		\end{equation}
		dos números reais (posto que ela traz o número $ \boldsymbol{0} $ como o ponto médio entre os números $ \boldsymbol{Z} $ e $ - \boldsymbol{Z} $), é justamente esse reconhecimento que acaba ratificando a interpretação dessa justaposição, da régua com a sua imagem que aparece na Figura \ref{regua-e-reta}\textbf{(a)}, como uma reta numérica.
		
		Já um outro ponto que reforça que o uso de espelhos planos é bem-vindo (dentro desse contexto de justificar todas as regras de sinais) segue de um comportamento que já foi explorado nas Figuras \ref{soma-total} e \ref{subtracao-total}. Afinal, sempre que um objeto é movimentado diante da superfície espelhada de tais espelhos, a sua imagem se move no sentido \textbf{oposto} ao longo da linha que une esse objeto à sua imagem. E observar esse comportamento é excelente posto que isso já consegue ilustrar, por exemplo, porque os números $ \boldsymbol{Z} $ e $ - \boldsymbol{Z} $ são interpretados como opostos um do outro.
		
		Aliás, note que essa oposição que existe $ \boldsymbol{Z} $ e $ - \boldsymbol{Z} $ já pode ser explorada para justificar, por exemplo, a regra de sinal que diz que \textquotedblleft mais com menos é menos\textquotedblright . Afinal, como o espelho plano está sobre o número $ \boldsymbol{0} $, fica bem fácil incutir a ideia de que $ \boldsymbol{0} $ é o ponto médio que existe entre esses dois números. Assim, como o resultado da média simples
		\begin{equation*}
			\frac{\boldsymbol{Z} + \left( - \boldsymbol{Z} \right) }{2}
		\end{equation*}
		precisa resultar nesse $ \boldsymbol{0} $, isso pode ser explorado para fazer com que os estudantes entendam tal regra de sinal dado que isso só ocorre se
		\begin{equation*}
			+ \left( - \boldsymbol{Z} \right) = - \boldsymbol{Z} \ .
		\end{equation*}
		De qualquer forma, é válido reforçar que pedir para que os estudantes observem que um objeto e a sua imagem se movem em sentidos opostos parece ser o mais lúdico para o bom entendimento dessas regras de sinais. Afinal, conforme já foi mencionado na Seção \ref{intro}, como todos os deslocamentos que caracterizam as somas de números negativos têm o mesmo sentido que os que caracterizam a subtração de números positivos, não é difícil usar isso para mostrar a validade da regra
		\begin{equation*}
			+ \left( - \boldsymbol{Z} \right) = - \left( + \boldsymbol{Z} \right) = - \boldsymbol{Z} \ .
		\end{equation*}
		Nestes termos, como esse último comentário pode ser perfeitamente estendido para mostrar a validade da regra
		\begin{equation*}
			- \left( - \boldsymbol{Z} \right) = + \left( + \boldsymbol{Z} \right) = + \boldsymbol{Z}
		\end{equation*}
		(já que os deslocamentos que caracterizam as somas de números positivos e as subtrações de números negativos têm a mesma direção), uma boa atividade que pode ser desenvolvida para fazer com que os estudantes entendam efetivamente todas essas regras de sinais é:
		\begin{itemize}
			\item introduzí-los ao conceito de \textquotedblleft oposto\textquotedblright , explorando atividades que envolvem o uso de espelhos onde eles podem observar as diversas imagens de objetos que são formadas;
			\item instigar que os estudantes avaliem se a distância que existe entre um objeto e um espelho plano é a mesma que existe entre a imagem e esse espelho, pedindo para que eles coloquem alguma régua perpendicularmente ao espelho; e
			\item fazer com que os estudantes percebam que o ato de colocar uma régua perpendicularmente a esse espelho, com o seu $ \boldsymbol{0} $ grudado na superfície espelhada, faz surgir uma reta que pode ser interpretada como a reta que contém os números reais.
		\end{itemize}
		Feito isso, é interessante pedir para que os estudantes desenhem uma semirreta e a coloquem perpendicularmente à superfície espelhada de um espelho plano pelo seu ponto $ \boldsymbol{0} $, tal como já foi ilustrado na Figura \ref{regua-e-reta}. Diante de tal construção, a segunda parte da atividade consiste em:
		\begin{itemize}
			\item pedir para que os estudantes coloquem alguns objetos sobre essa semirreta (como borrachas, apontadores etc.) e, anotando as posições $ \boldsymbol{p_{\mathrm{im}}} $ e $ \boldsymbol{p_{\mathrm{ob}}} $ que esses objetos e as suas imagens assumem respectivamente, percebam que $ \boldsymbol{p_{\mathrm{im}}} = - \boldsymbol{p_{\mathrm{ob}}} $;
			\item pedir para que os estudantes movam esses objetos sobre essa semirreta, e percebam que os movimentos que estão associados às operações de somas/subtração sempre afastam/aproximam esses objetos e as suas imagens do espelho;
			\item anotem os sentidos dos movimentos que estão asssociados às
			\begin{itemize}
				\item[\textbf{(a)}] somas de números positivos,
				\item[\textbf{(b)}] subtrações de números negativos,
				\item[\textbf{(c)}] somas de números negativos e
				\item[\textbf{(d)}] subtrações de números positivos
			\end{itemize}
			e identifiquem quais deles têm os mesmos sentidos.
		\end{itemize}
		Como consequência da identificação que é pedida neste último item, todas essas regras de sinais poderão ser justificadas.
		
		\subsection{Uma observação física}
		
			Entretanto, vale fazer uma observação bem importante aqui, a qual tem a ver com uma coisa que não costuma ser muito bem explorada no Ensino Básico: \textbf{a Matemática é a linguagem pela qual a Física se escreve}. A propósito, apesar de eu saber, por exemplo, que para enviar um foguete ou um satélite para o espaço exigia saber fazer muitas contas, eu mesmo, quando eu estava no primeiro ano do Ensino Médio, ainda não entendia muito bem por que o meu professor de Matemática sabia tantas coisas sobre Física. Afinal, ele era apenas um professor de Matemática e aquela coisa que chamavam de Física, e que me foi apresentada até o final do meu Ensino Médio, não era a Ciência que eu também já amava: ela era só mais uma matéria da escola, onde todos os seus professores apenas me pediam para eu decorar um monte de fórmulas sem nunca me dar qualquer explicação sobre o porquê de todas elas.
			
			E por que fazer essa observação, sobre a Matemática ser a linguagem pela qual a Física se escreve, é importante? Porque, quando nós voltamos as nossas atenções exatamente para o contexto deste artigo, não é errado afirmar que o fato do eixo principal de um espelho esférico corresponder à reta dos números reais não se deve a uma mera convenção: essa correspondência é \textbf{proposital} posto que a Matemática, por ser a linguagem pela qual a Física se escreve, mantém uma correspondência natural para com a realidade, para com tudo aquilo que nos cerca. Ou seja, apesar de algumas pessoas acharem que a Matemática é apenas uma concepção do cérebro humano, o fato é que ela pode mesmo ser observada em diversas situações na nossa realidade e uma das situações onde isso fica bem claro é, curiosamente, diante de um espelho. Afinal, se assim não fosse, nós jamais enxergaríamos a reta dos números reais ao colocar uma simples régua na frente de um espelho plano.
			
			Nestes termos, uma atividade que pode ser realizada em sala de aula, não apenas pelos professores de Matemática do Ensino Básico, mas (principalmente) pelos professores de Física nas aulas que abordam a formação de imagens por espelhos esféricos, é exatamente essa que está sendo proposta neste artigo. É claro que, numa aula de Física que versa sobre essa formação de imagens, as principais preocupações são outras -- como explorar as características que as imagens possuem, se elas são maiores ou menores, invertidas ou direitas etc. Porém, realizar uma atividade que também tem um caráter um pouco mais matemático numa aula de Física, para que os estudantes notem (ao menos) que a reta dos números reais pode ser construída com o auxílio de uma régua e de um espelho plano, é algo muito bem-vindo. Afinal de contas, além dessa atividade ser bastante útil para mostrar (e explorar) a presença da Matemática de um jeito que os estudantes não esperam, é justamente essa presença que acaba deixando claro por que, por exemplo, todas as posições atrás dos espelhos são descritas apenas por números negativos haja vista que essa é uma dúvida bastante comum entre os estudantes.
			
	\section{\label{conclusions}Comentários finais}
	
		Embora a história que eu contei ao longo deste artigo traga uma experiência bastante pessoal (onde eu também exponho, por exemplo, a minha opinião sobre eu achar muito chatas boa parte das aulas de Física e de Matemática que eu tive no Ensino Básico), o fato é que essa realidade do \textquotedblleft decorar\textquotedblright \hspace*{0.01cm} e do \textquotedblleft reproduzir\textquotedblright \hspace*{0.01cm} (que era justamente o que eu achava muito chato) ainda é bastante recorrente nas aulas da grande maioria das escolas de Ensino Básico do Brasil. Entretanto, na grande maioria das vezes, essa realidade não chega nem a ser culpa direta dos professores que lecionam todas essas disciplinas, mas, sim, do sistema educacional onde todos estão imersos, o qual não é capaz de atrair bons professores e de, tão pouco, valorizar os que já estão dando aulas oferencendo, a eles, uma melhor formação.
		
		Um bom exemplo disso é o fato de que uma parte considerável dos professores que me deram aulas no Ensino Básico sequer tinham a formação adequada para lecionar todas as matérias que eles me lecionaram. Afinal de contas, eu já tive aulas de Biologia com um farmacêutico, de Química com um segundo anista de Medicina e quase todos os professores que me deram aulas de Matemática nas (antigas) quinta, sexta e sétima séries do Ensino Básico eram formados em Ciências Biológicas e não em Matemática. Aliás, essa também era a formação da professora que me deu aulas de Física no primeiro ano do Ensino Médio: ela era formada em Ciências Biológicas, não em Física. E, quando eu perguntei para ela, em sala de aula, sobre o porquê de existir um fator $ \boldsymbol{1} / \boldsymbol{2} $ na equação horária
		\begin{equation*}
			\boldsymbol{s} \left( \boldsymbol{t} \right) = \boldsymbol{s_{0}} + \boldsymbol{v_{0} t} + \frac{\boldsymbol{1}}{\boldsymbol{2}} \ \boldsymbol{at^{2}}
		\end{equation*}
		do movimento uniformemente acelerado, ela simplesmente \textbf{riu} (muito provavelmente para disfarçar que ela não sabia a resposta), mas não foi só isso: além de rir, ela disse que \textbf{não era para eu me preocupar com isso pois nunca, na minha vida, eu ia precisar saber por que esse fator $ \boldsymbol{1} / \boldsymbol{2} $ estava ali}. Alguns anos depois, eu obtive o meu \textbf{doutorado em Física}.
		
		É claro que, às vezes, eu até chego a ter vontade de reencontrar essa professora só para lembrá-la deste episódio e poder explicar para ela o porquê daquele fator $ \boldsymbol{1} / \boldsymbol{2} $. Afinal, além de ser bastante possível que ela ainda não saiba a resposta da pergunta que eu lhe fiz, eu subentendo que, por trás do que ela me disse, talvez existisse um certo preconceito já que, como boa parte dos estudantes daquela escola eram filhos de empregadas domésticas, pedreiros e de pais com outras profissões que (infelizmente) são vistas como menores pela atual sociedade brasileira (como também era o meu caso), provavelmente ela supunha que ninguém, naquela sala de aula, iria tão longe nos estudos. Ou seja, por mais que a maioria dos professores realmente não tenha culpa alguma das eventuais falhas que existem nas suas formações, não existe motivo algum para que eles duvidem da capacidade que os estudantes têm para aprender alguma coisa, muito menos se essa dúvida eventualmente estiver embebida em algum preconceito.
		
		Aliás, duvidar da capacidade que os estudantes possuem, por qualquer que seja o motivo, ou mesmo desincentivar qualquer questionamento pertinente que eles façam dentro de uma sala de aula do Ensino Básico, chega até mesmo a ser contraproducente à luz do que, por exemplo, afirma a \textbf{teoria do desenvolvimento cognitivo} de Jean Piaget \cite{piaget1,piaget2,piaget3}. Afinal de contas, é justamente na época que os estudantes adentram na segunda metade do Ensino Fundamental que, segundo essa teoria, a maioria deles já começa a ser capaz de solucionar problemas criando conceitos e ideias, assim como se valendo de pensamentos formais abstratos. Ou seja, por mais que existam respostas que exijam, de fato, um conhecimento um pouco mais profundo do que o que é oferecido no Ensino Básico, muito provavelmente esses estudantes já são perfeitamente capazes não apenas de entender todas as respostas que lhes sejam apresentadas de um jeito intelígivel, mas de, inclusive, obtê-las por conta própria desde que uma discussão seja fomentada para isso (note que a história real que eu acabo de contar neste texto é um bom exemplo disso). E é justamente dentro desse contexto, de fomentar uma discussão, que esse artigo, que você, leitor, está acabando de ler agora, se encaixa. Afinal, basta ver que, além dele trazer uma proposta bastante lúdica através da qual um professor pode abordar esse tema das regras de sinais, essa mesma proposta também pode fomentar diversas discussões que até podem fugir do escopo desse tema, mas que ainda estão dentro de um nicho que pertence à Física e à Matemática. Aliás, uma das outras discussões que o uso de espelhos no ensino da Matemática pode fomentar versa, por exemplo, sobre o conceito matemático de \textbf{conjunto imagem} que está associado a uma função, já que não é à toa que esse conjunto possui justamente o nome daquilo que um espelho é capaz de formar. Porém, como qualquer comentário adicional que eu possa fazer sobre isso só vai estender ainda mais esse texto (que já está um pouco longo) sem a menor necessidade, eu prefiro deixar para falar sobre isso num artigo futuro.
			
	\section{Agradecimentos}
	
		Este trabalho infelizmente não foi financiado por nenhum órgão, seja ele privado ou governamental, pois, ao contrário do que acontece em algumas das grandes nações, ser cientista e professor aqui no Brasil ainda parece ser um ato de heroísmo e de afronta ao sistema que está aqui estabelecido há séculos.
		
		No entanto, eu agradeço profundamente a todos os personagens que fizeram parte, seja direta ou indiretamente, desse pedaço da minha história que eu apresentei neste artigo. Afinal, se não fosse por causa de todos eles, e de todos os desafios felizes e não felizes que eles colocaram na minha frente, eu jamais teria obtido todas as provas de que eu realmente estava no caminho certo, rumo ao \textquotedblleft saber a fazer todas as contas possíveis, principalmente todas aquelas que permitiam entender as coisas do Universo\textquotedblright , mesmo que esse sistema já conspirasse, desde aquela época, para que outras pessoas como eu não conseguissem ter acesso a tais informações. Este jamais seria o meu caso pois, além de eu sempre ter acreditado na minha capacidade (e nunca duvidado dela), outra pessoa também sempre acreditou: a minha mãe, Celsina Jacinta de Araujo, a quem eu agradeço profundamente não apenas por toda a ajuda na construção da minha trajetória, mas por ter (inclusive) saído em minha defesa, protagonizando uma das principais cenas que eu contei neste artigo.
		
		De qualquer forma, apesar de eu não ter citado os nomes dos demais personagens, eu preciso citar o nome de duas pessoas que, de alguma maneira, participaram (in)diretamente da construção deste artigo. E uma dessas pessoas é o Prof. Dr. Vinicio de Macedo Santos (FE-USP) haja vista que foi durante uma das suas aulas de \textquotedblleft EDM0427 -- Metodologia do Ensino da Matemática I\textquotedblright \hspace*{0.01cm} que eu, há exatos dezesseis anos ainda como estudante de graduação, apresentei (pela primeira vez) essa minha abordagem para o ensino das regras de sinais. Já a outra pessoa que eu preciso agradecer aqui é o meu grande amigo, o Prof. Dr. Leandro Daros Gama (IFSP), pela leitura minuciosa deste artigo, assim como pelas adicionais sugestões que fizeram desse texto algo bem melhor de ser lido.
		
		Todavia, antes de eu terminar este texto, não é muito justo eu deixar de agradecer algumas pessoas que, apesar de eu sequer conhecê-las, disponibilizaram algumas imagens extremamente bonitas que fizeram parte da apresentação deste artigo. E as primeiras delas são as que, usando os codinomes \textquotedblleft photosforyou\textquotedblright , \textquotedblleft fireboltbyl\textquotedblright \hspace*{0.01cm} e \textquotedblleft S. Hermann \& F. Richter\textquotedblright , disponibilizaram publicamente as fotografias que constam nas Figuras \ref{mirror-balls}, \ref{spherical-mirror} e \ref{jaguar} respectivamente \cite{photosforyou,fireboltbyl,pixel2013}. Já as outras pessoas que também merecem o meu agradecimento aqui são as que, com os codinomes \textquotedblleft JJuni\textquotedblright , \textquotedblleft OpenClipart-Vectors\textquotedblright \hspace*{0.01cm} e \textquotedblleft Clker-Free-Vector-Images\textquotedblright , também disponibilizaram publicamente algumas imagens que eu acabei editando para criar as Figuras \ref{regua-e-imagem}, \ref{regua-e-reta}, \ref{soma-total}, \ref{subtracao-total}, \ref{convencao} e \ref{imagem-real} com o auxílio de uma programação em \LaTeX \hspace*{0.01cm} \cite{jjuni,openclip-1,openclip-2,clker-free}. E, já que eu acabei de falar em \LaTeX , um último agradecimento também cabe a quem, usando o codinome \textquotedblleft user121799\textquotedblright , publicou o código  \textquotedblleft .tex\textquotedblright \hspace*{0.01cm} que consta na Ref. \cite{user121799}, o qual eu usei para gerar a Figura \ref{corte-esfera}.
	
	\appendix
	
	\section{\label{appendix}Sobre a formação de imagens por um espelho esférico}
	
		Se existe uma coisa, hoje, que já está muito bem estabelecida é que a luz pode ser interpretada como uma \textbf{onda eletromagnética}, posto que, a ela, estão associados dois campos (um elétrico e outro magnético) que satisfazem uma equação de onda \cite{moyses}. É claro que o objetivo aqui, com este Apêndice, não é o de pormenorizar qualquer um dos motivos que são capazes de justificar essa interpretação ondulatória da luz. No entanto, não deixa de ser interessante destacar (mesmo que pareça ser apenas a título de uma curiosidade) que levou muito tempo para convencer toda a comunidade científica de que a luz realmente era uma onda. No caso, todos os Séculos XVI e XVII, assim como mais da metade do Século XVIII, foram voltados para uma discussão intensa entre aqueles que defendiam essa interpretação, entre os quais podemos destacar os nomes de Christiaan Huygens (1629 -- 1695), Thomas Young (1773 -- 1829) e Augustin Jean Fresnel (1788 -- 1827) \cite{crew}, e aqueles que defendiam uma outra interpretação para a luz: a de que ela poderia ser interpretada \textbf{corpuscularmente}. E entre os que defendiam essa última interpretação, o grande nome era o de Isaac Newton (1643 -- 1727) \cite{newton-optica}. 
		
		Aliás, o fato de Newton ser partidário dessa interpretação corpuscular para a luz foi, certamente, a maior dificuldade que precisou ser superada por todos os que defendiam a interpretação ondulatória. Afinal de contas, como ir contra o que um dos nomes mais influentes e respeitados da história da Física (já em vida) afirmava? Realmente foi uma tarefa muito difícil, mesmo com Young tendo desenvolvido, já no início do Século XIX (isto é, mais de meio século depois que Newton havia falecido), o seu \textbf{experimento da fenda dupla} \cite{young-light-1,young-light-2} onde era bem claro o fenômeno de \textbf{interferência da luz}, um fenômeno que só podia ser justificado se a luz fosse uma onda \cite{young}. Ou seja, essa aceitação, de que a luz realmente poderia ser interpretada como uma onda, foi extremamente lenta e só foi finalmente ratificada \textbf{(i)} com a ajuda de James Clerk Maxwell (1831 -- 1879) que, com os seus estudos físicos e matemáticos sobre a natureza da Eletricidade e do Eletromagnetismo, concluiu, em 1864, que a luz era uma forma de onda eletromagnética \cite{maxwell-wave} e \textbf{(ii)} com a verificação experimental das equações que Maxwell obteve nos seus estudos, a qual só foi feita em 1887 por Heinrich Rudolf Hertz (1857 -- 1894) \cite{hertz}.
		
		Entretanto, vale notar que Newton defendia uma interpretação cospuscular não porque ele queria ser \textquotedblleft do contra\textquotedblright . Newton era um cientísta excepcional e, se ele defendia tal interpretação, era porque, além da interpretação ondulatória da luz possuir uma falha bastante estrutural na sua época, também existiam diversos indícios experimentais que pareciam endossar uma interpretação corpuscular\footnote{Essa falha \textquotedblleft bastante estrutural\textquotedblright \hspace*{0.01cm} se pauta no fato de que os partidários da interpretação ondulatória acreditavam que a luz era uma onda \textbf{longitudinal}, uma vez que todas as ondas conhecidas até então tinham essa propriedade. E como alguns estudos já deixavam bem claro que, se a luz fosse uma onda, ela não poderia ser longitudinal já que isso não se ajustava ao fenômeno da sua polarização, isso acabava sendo mais um dos argumentos usados pelos defensores da interpretação corpuscular. A conclusão de que a luz era uma onda \textbf{transversal} só foi apresentada, pela primeira vez, por Young em 1817: ele chegou a essa conclusão depois de analisar os trabalhos que Fresnel e Dominique François Arago (1786 -- 1853) realizaram, os quais afirmavam que dois feixes de luz com polarizações perpendiculares nunca se interferem \cite{ribeiro}. Tempos mais tarde foi Maxwell que, com o seu trabalho seminal de 1864, acabou demostrando que a transversalidade associada à luz se deve ao fato dos seus campos elétrico e magnético oscilarem \textbf{perpendicularmente} à direção por onde essa luz se propaga \cite{maxwell-wave}.}. E um desses indícios era, por exemplo, o fato de que, quando qualquer feixe de luz incide sobre uma superfície bem polida, ele sempre é refletido do mesmo jeito que uma partícula é refletida por uma parede lisa e rígida: ou seja, com o seu ângulo de reflexão tendo o mesmo valor que o de incidência conforme ilustra a Figura \ref{bolinha}. Assim, diante não apenas desse indício, mas de outros que também apontavam para uma correspondência entre a Óptica e a Mecânica, não havia como Newton não achar válida a ideia de que a luz podia ser interpretada em termos corpusculares, muito embora ele também não tenha deixado de reconhecer que alguns fenômenos ópticos exigiam uma interpretação ondulatória dentro de algum contexto \cite{salvetti}.
		
		\subsection{Como os feixes de luz são refletidos por um espelho esférico?}
		
			Aliás, todos esses indícios, que faziam com que Newton interpretasse a luz corpuscularmente, eram tão robustos que, até hoje, eles nos permitem ignorar completamente o fato da luz ser uma onda em diversas situações. E uma dessas situações é aquela em que nós precisamos entender como a imagem de um objeto é formada por um espelho esférico.
			\begin{figure}[!t]
				\centering
				\includegraphics[viewport=330 10 0 155,scale=1.3]{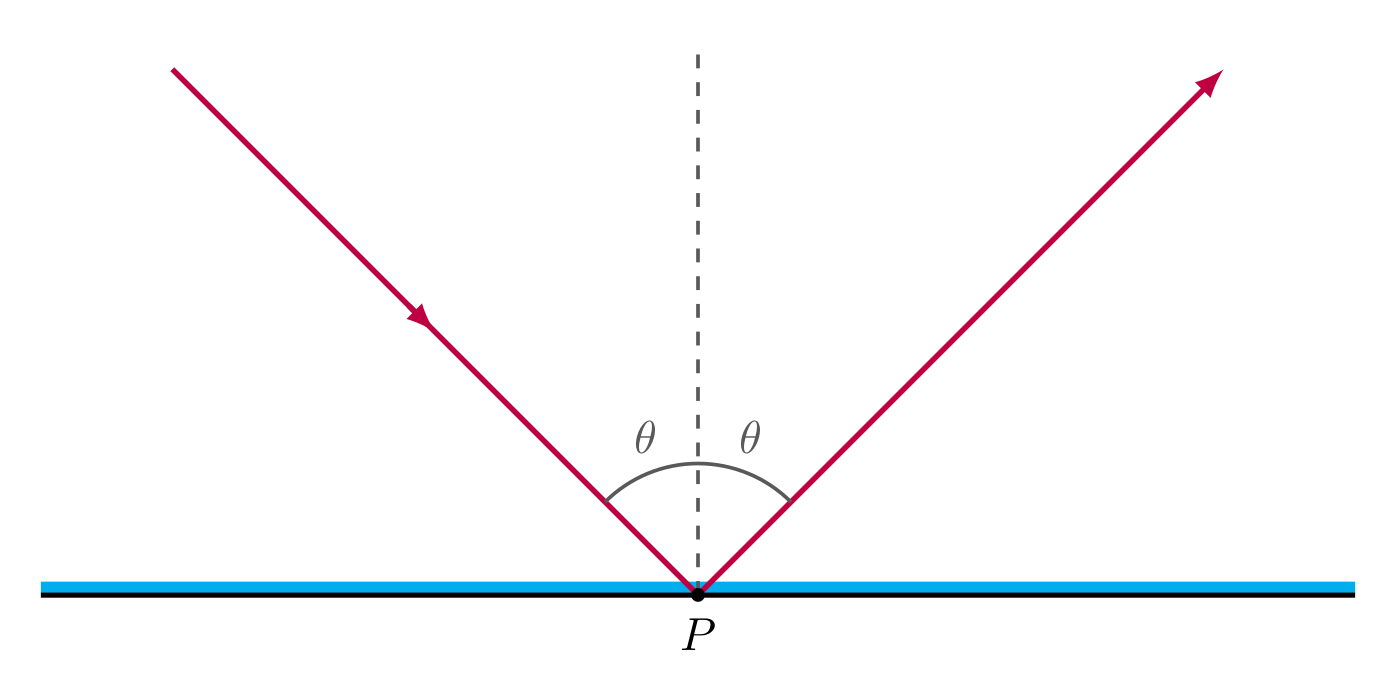}
				\caption{\label{bolinha}Situação de um feixe de luz quando refletido por uma superfície espelhada num ponto $ \boldsymbol{P} $. Note que o feixe incidente (aquele que está à esquerda da linha hachurada) forma um ângulo $ \boldsymbol{\theta } $ com uma linha que é normal ao espelho (no caso, a mesma linha hachurada), ângulo o qual possui o mesmo valor que aquele que o feixe refletido (que aparece à direita) forma com essa mesma normal.}
			\end{figure}
			
			Só que antes de entender como isso acontece e de, consequentemente, responder à questão que dá nome a esta Subseção, é preciso deixar claro um fato muito importante que versa sobre a capacidade que nós, humanos, temos de enxergar um objeto com os nossos olhos. E que fato é esse? É o fato de que, além de ser necessário que os nossos olhos estejam em bom funcionamento para isso, também é necessário que haja \textbf{luz} sobre esse objeto. Afinal de contas, nós só enxergamos um objeto porque, depois que os raios da luz que o ilumina são (por ele) refletidos, uma parte desses raios adentra os nossos olhos e estimulam algumas células que ficam nas nossas retinas. E como esses estímulos, que ocorrem sob a forma de impulsos elétricos, acabam sendo transmitidos do olho para o nosso cérebro através dos nervos ópticos, é possível afirmar que o ato de enxergar é um ato não apenas de detectar luz, mas de traduzir toda a informação desse conjunto de impulsos elétricos para o modo gráfico. É exatamente esse gráfico que nós enxergamos e não o que realmente existe ao nosso redor, uma vez que o espectro do comprimento das ondas eletromagnéticas que nós, humanos, somos capazes de enxergar é bastante diminuto \cite{okuno}.
			
\begin{figure}[!t]
				\centering
				\includegraphics[viewport=330 10 0 160,scale=1.3]{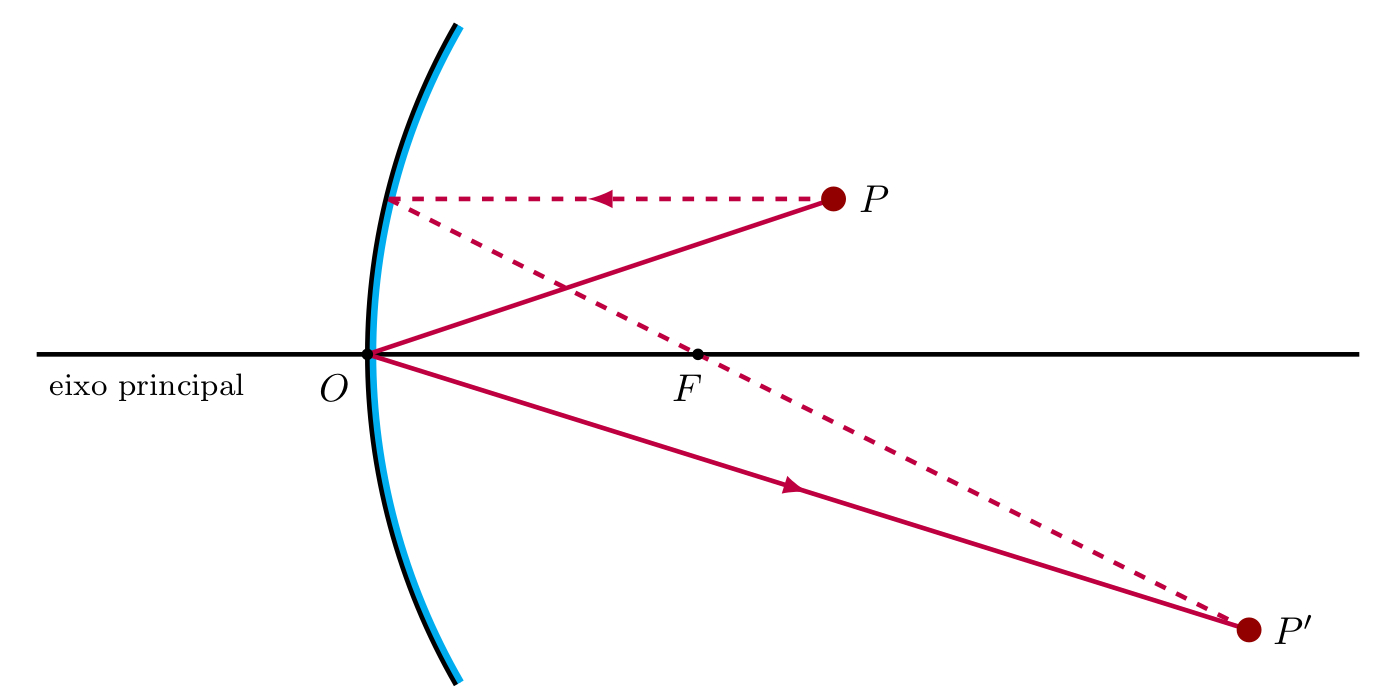}
				\caption{\label{esquema}Esquema que é tradicionalmente utilizado para obter a imagem de um ponto luminoso $ \boldsymbol{P} $ que está diante de um espelho côncavo. Aqui, estão destacados apenas dois raios de luz que saem desse ponto: um (que está destacado por uma linha contínua na cor púrpura) que incide sobre o ponto $ \boldsymbol{O} $ por onde passa o eixo principal; e outro (que é destacado por uma linha hachurada com essa mesma cor púrpura) que, por ser paraxial, está sendo refletido para o foco do espelho. No caso da imagem $ \boldsymbol{P^{\prime }} $ desse ponto luminoso, ela é identificada como o ponto de cruzamento desses dois raios de luz depois deles serem refletidos por esse espelho.}
			\end{figure}

			Outra coisa que também é importante de ser destacada aqui é que, quando um objeto está sendo iluminado, cada um dos seus pontos refletem raios de luz nas mais diversas direções. Nestes termos, é bastante plausível reconhecer que, quando um objeto iluminado é colocado diante de um espelho esférico, sempre existe um raio de luz que sai de um dos pontos desse objeto e chega até o ponto de intersecção desse espelho com o seu eixo principal. Essa é justamente a situação aparece na Figura \ref{esquema},
			onde, devido ao fato do eixo principal ser perpendicular ao espelho, é possível notar que esse mesmo raio de luz está sendo refletido nos mesmos moldes já observados na Figura \ref{bolinha}: ou seja, com o ângulo de reflexão tendo o mesmo valor que o de incidência.
			
\begin{figure}[!t]
				\centering
				\includegraphics[viewport=330 10 0 180,scale=1.3]{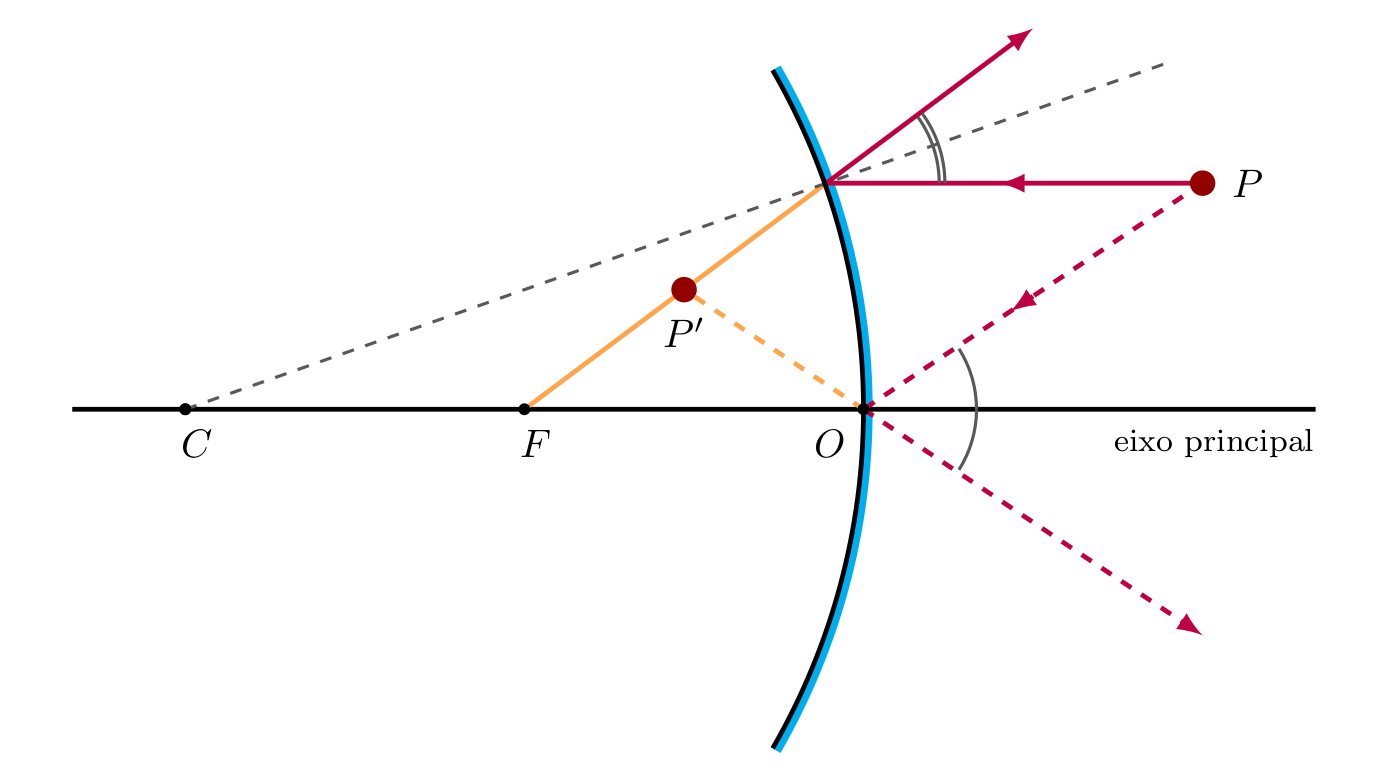}
				\caption{\label{extensao}Esquema tradicionalmente usado para obter a imagem de um ponto luminoso $ \boldsymbol{P} $ que está diante de um espelho \textbf{convexo} (ou seja, de um espelho cuja face espelhada (destacada em ciano) é a superfície externa da calota que lhe deu origem). Aqui, novamente estão destacados apenas dois raios de luz: um (que está destacado por uma linha hachurada com a cor púrpura) que sai de $ \boldsymbol{P} $ e incide sobre o ponto $ \boldsymbol{O} $ por onde passa o eixo principal; e outro paraxial (sendo destacado por uma linha contínua com a mesma cor púrpura) que também sai de $ \boldsymbol{P} $ e está sendo refletido para o nordeste que quem lê esta página. Note que os dois raios luminosos refletidos \textbf{não se cruzam}. No entanto, como a reta (hachurada em cinza) que passa pelo ponto $ \boldsymbol{C} $ pode ser interpretada como uma \textbf{bissetriz} para o raio de luz paraxial e o que é refletido para o nordeste, a extensão deste raio refletido (que está destacada com a cor laranja) cruza o eixo principal exatamente no seu ponto $ \boldsymbol{F} $. Nestes termos, como a extensão do primeiro raio refletido (o que foi destacado por uma linha contínua na cor púrpura) intercepta a primeira extensão num ponto $ \boldsymbol{P^{\prime }} $, a mesma lógica que foi empregada na figura anterior nos permite afirmar que é justamente esse $ \boldsymbol{P^{\prime }} $, que está atrás do espelho, que pode ser interpretado como a imagem de $ \boldsymbol{P} $.}
			\end{figure}

			Todavia, como cada ponto de um objeto iluminado reflete raios de luz nas mais diversas direções, um outro raio que também pode ser destacado aqui é aquele se propaga rumo ao espelho \textbf{paralelamente} ao seu eixo principal. E, conforme já foi observado na Subseção \ref{subsection:paraxiais}, se esse raio for \textbf{paraxial}, ele será refletido do espelho diretamente para o seu foco. Desta maneira, uma boa maneira de reconhecer onde a imagem deste ponto será construída é avaliando, por exemplo, onde esses dois raios se cruzam. No caso que aparece na própria Figura \ref{esquema}, esse cruzamento se dá na frente da superfície espelhada e, por isso, a imagem construída é reconhecida como \textbf{real}. Já no caso da situação que aparece na Figura \ref{extensao}, esse cruzamento se dá atrás do espelho e, por isso, a imagem é reconhecida como \textbf{virtual}.
			
\begin{figure}[!t]
				\centering
				\includegraphics[viewport=460 10 0 295,scale=1.0]{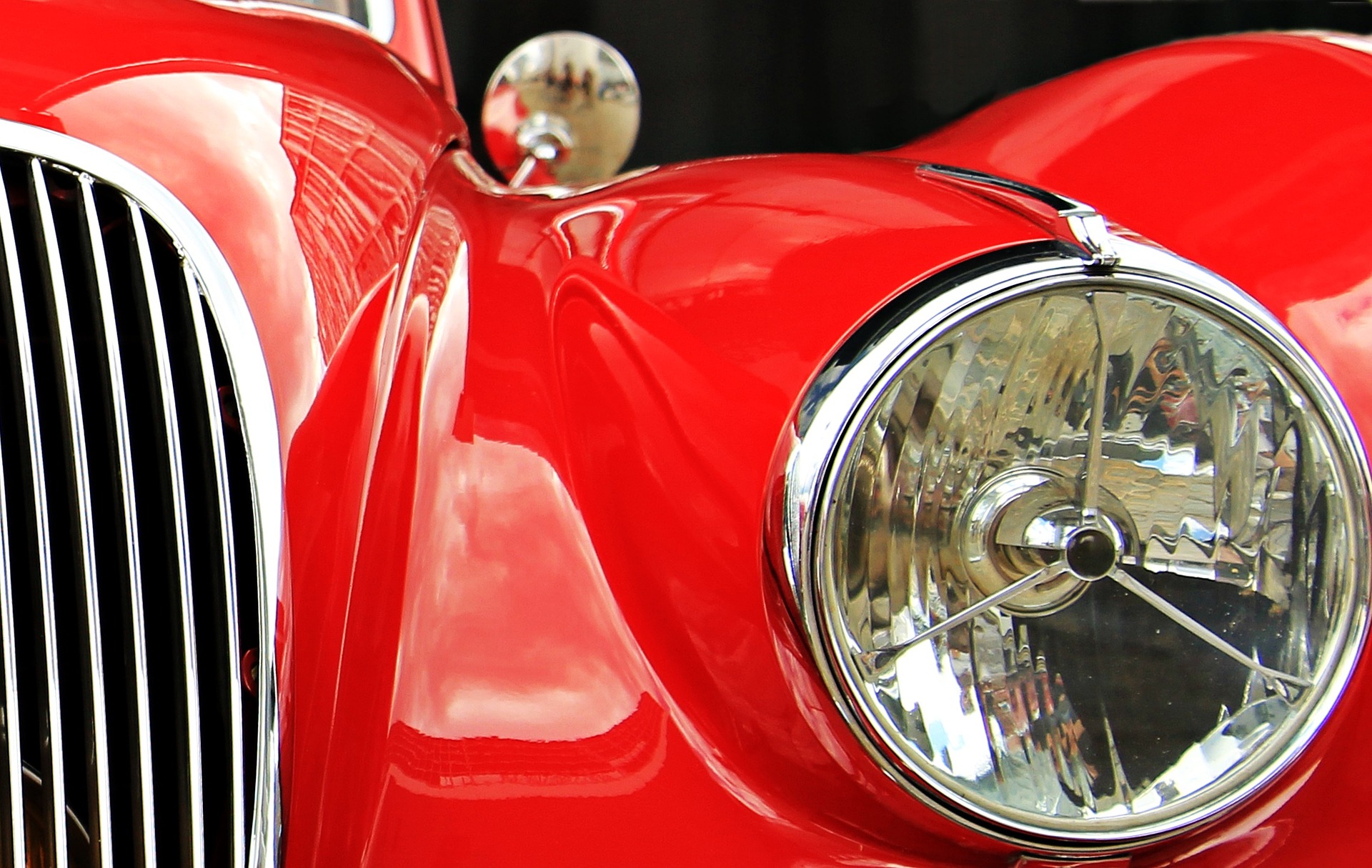}
				\caption{\label{jaguar}Fotografia de um antigo Jaguar vermelho onde é dado destaque a um dos seus faróis. Note que, no interior desse farol, constam um espelho côncavo e uma fonte de luz que está sendo sustentada por três hastes. Essa fonte de luz está exatamente no foco do espelho côncavo para maximizar o poder de iluminação do farol.}
			\end{figure}

			Note que o caso da imagem que aparece na Figura \ref{light-reflection} também pode ser justificado nos mesmos moldes que esses. Todavia, é importante notar que a identificação da imagem na Figura \ref{light-reflection} foi feita com a ajuda de um terceiro exemplo de raio (o que está destacado em hachurado nessa mesma Figura \ref{light-reflection}), que é aquele que sai do objeto iluminado e que passa diretamente pelo foco do espelho antes de chegar à superfície desse espelho. Ou seja, esse terceiro exemplo de raio luminoso está fazendo um caminho que é \textquotedblleft inverso\textquotedblright \hspace*{0.01cm} ao dos raios paraxiais que são refletidos por um espelho para o seu foco. Só que, por mais estranho que possa parecer identificar uma imagem destacando esse terceiro exemplo de raio luminoso, vale notar que a existência desse raio já é muito bem conhecida, por exemplo, por aqueles que trabalham com luz, especialmente por aqueles que constroem alguns dispositivos de iluminação, entre os quais podemos destacar muitas lanternas, faróis de carro (principalmente os mais antigos) e, principalmente, os holofotes. Afinal de contas (e conforme bem ilustra a Figura \ref{jaguar} a seguir), a lógica por trás da construção desses dispositivos se pauta na alocação de uma fonte de luz (que possui o menor tamanho possível) exatamente no foco de um espelho côncavo uma vez que, como todos os raios dessa luz (que incidem sobre a superfície espelhada) são refletidos para longe do espelho, é justamente isso que maximiza a iluminação.
		
		\subsection{Obtendo a equação de Gauss}
		
			Diante de tudo o que acabou de dito, sobre esses raios luminosos que saem de um objeto iluminado e rumam para um espelho esférico, fica mais fácil demonstrar, por exemplo, a validade da equação de Gauss
			\begin{equation*}
				\frac{\boldsymbol{1}}{\boldsymbol{p_{\mathrm{ob}}}} + \frac{\boldsymbol{1}}{\boldsymbol{p_{\mathrm{im}}}} = \frac{\boldsymbol{1}}{\boldsymbol{f}} \ .
			\end{equation*}
			E lembrando que $ \boldsymbol{p_{\mathrm{ob}}} $ e $ \boldsymbol{p_{\mathrm{im}}} $ denotam respectivamente as posições que o objeto $ \boldsymbol{P} $ e a sua imagem $ \boldsymbol{P^{\prime }} $ assumem sobre o eixo principal de um espelho, é possível afirmar que já existem duas coisas, que constam na Figura \ref{esquema}, que nos ajudam a demonstrar tal equação\footnote{Embora a Figura \ref{esquema} tenha sido escolhida para nos guiar nesta demostração, vale notar que qualquer outra figura, que traga um espelho esférico e qualquer ponto luminoso com a sua imagem (identificada através do cruzamento de dois raios refletidos que também constem nessa figura), pode ser usada para este fim. No caso da escolha que foi feita pela Figura \ref{esquema}, ela só foi feita porque, além dessa figura já ter sido construída previamente, ela é bem menos poluída visualmente do que, por exemplo, a Figura \ref{extensao}.}.
			
\begin{figure}[!t]
				\centering
				\includegraphics[viewport=330 10 0 165,scale=1.3]{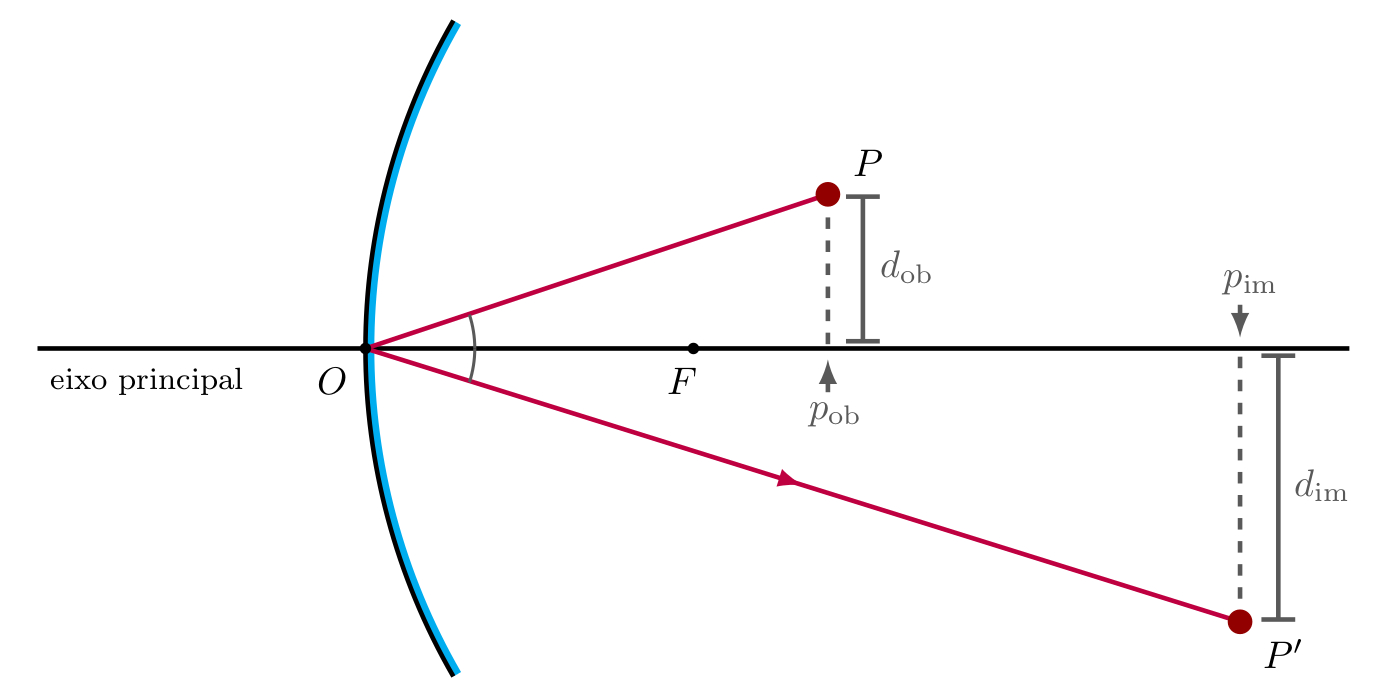}
				\caption{\label{triangulos-1}Situação já exposta na Figura \ref{esquema} onde, agora, é dada ênfase ao fato dos ângulos, que estão subtendidos pelo eixo principal e os raios luminosos incidente (que está acima do eixo) e refletido (que está abaixo do eixo), terem o mesmo valor. Afinal de contas, como esse fato garante a igualdade entre as tangentes desses ângulos (ou seja, entre as razões dos catetos oposto e adjacente dos dois triângulos que aparecem nesta figura), é exatamente isso que justifica a relação (\ref{altura-1}).}
			\end{figure}

\begin{figure}[!t]
				\centering
				\includegraphics[viewport=330 10 0 165,scale=1.3]{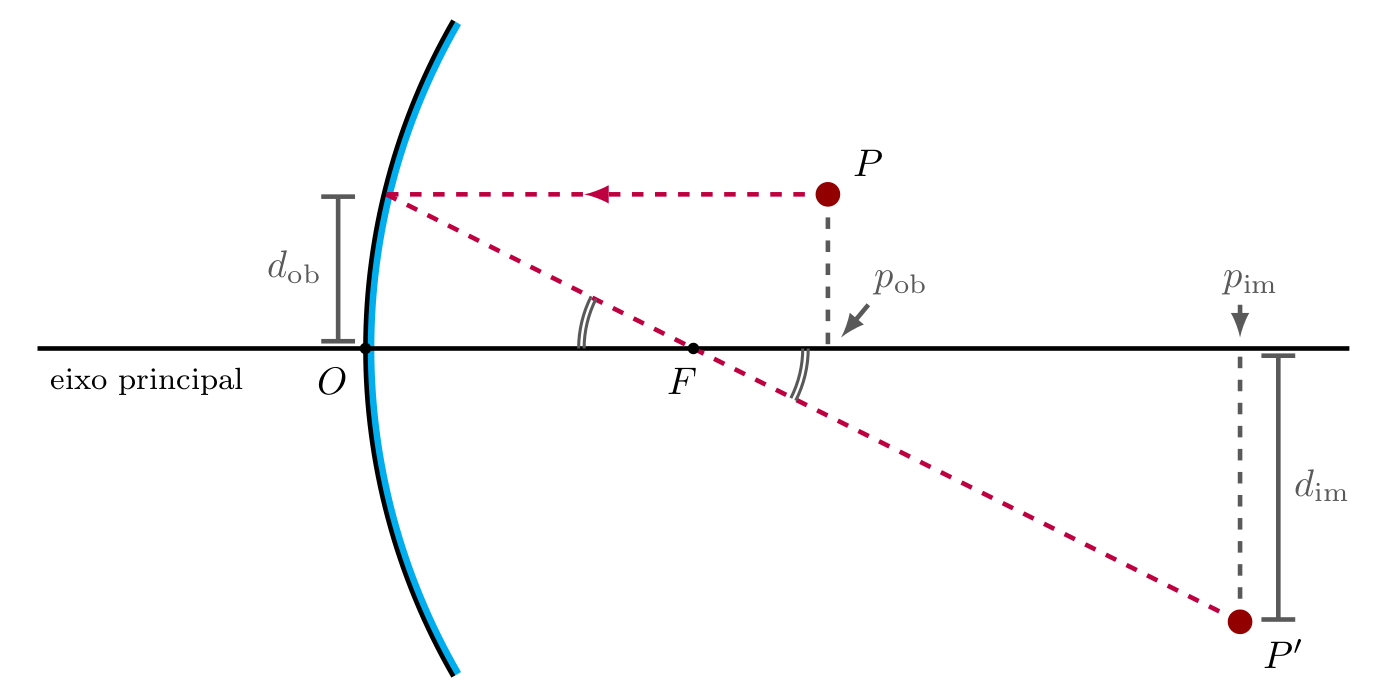}
				\caption{\label{triangulos-2}Novamente estamos diante da mesma situação que já foi exposta na Figura \ref{esquema} só que, desta vez, dando destaque ao fato do raio luminoso refletido interceptar o eixo principal no seu ponto $ \boldsymbol{F} $. Afinal, como os ângulos que ficam subtendidos por esse raio e pelo eixo têm o mesmo valor, é justamente isso que leva à relação (\ref{altura-1}) desde que o raio luminoso incidente seja paraxial. Note que, como, no triângulo à esquerda de $ \boldsymbol{F} $, o cateto adjacente mede $ \boldsymbol{f} $, o cateto que é adjacente no triângulo à direita mede $ \boldsymbol{p_{\mathrm{im}}} - \boldsymbol{f} $.}
			\end{figure}

			A primeira delas é que, se nós considerarmos que $ \boldsymbol{d_{\mathrm{ob}}} $ e $ \boldsymbol{d_{\mathrm{im}}} $ são as respectivas distâncias que o objeto e a sua imagem mantém do eixo principal, passa a ser válido afirmar que
			\begin{equation}
				\frac{\boldsymbol{d_{\mathrm{ob}}}}{\boldsymbol{d_{\mathrm{im}}}} = \frac{\boldsymbol{p_{\mathrm{ob}}}}{\boldsymbol{p_{\mathrm{im}}}} \label{altura-1}
			\end{equation}
			haja vista que o ângulo de incidência e de reflexão, do raio que está sendo destacado com uma linha contínua, são \textbf{idênticos}. A Figura \ref{triangulos-1} nos ajuda a entender um pouco melhor o porquê desta relação.
			Já se nós olharmos para o raio paraxial e para a sua reflexão (ambos hachurados nessa mesma Figura \ref{esquema}), também não é difícil perceber que		
			\begin{equation}
				\frac{\boldsymbol{d_{\mathrm{ob}}}}{\boldsymbol{d_{\mathrm{im}}}} = \frac{\boldsymbol{f}}{\boldsymbol{p_{\mathrm{im}} - \boldsymbol{f}}} \ , \label{altura-2}
			\end{equation}
			uma vez que esse raio intercepta o eixo principal no seu ponto $ \boldsymbol{F} $ (vide a Figura \ref{triangulos-2} para entender como essa relação fica justificada).
			Nestes termos, como as relações (\ref{altura-1}) e (\ref{altura-2}) nos mostram que			
			\begin{equation*}
				\frac{\boldsymbol{p_{\mathrm{ob}}}}{\boldsymbol{p_{\mathrm{im}}}} = \frac{\boldsymbol{f}}{\boldsymbol{p_{\mathrm{im}} - \boldsymbol{f}}} \ ,
			\end{equation*}
			uma simples manipulação algébrica nos leva a			
			\begin{equation*}
				\boldsymbol{p_{\mathrm{im}}} \cdot \boldsymbol{f} = \boldsymbol{p_{\mathrm{ob}}} \left( \boldsymbol{p_{\mathrm{im}}} - \boldsymbol{f} \right) = \boldsymbol{p_{\mathrm{ob}}} \cdot \boldsymbol{p_{\mathrm{im}}} - \boldsymbol{p_{\mathrm{ob}}} \cdot \boldsymbol{f}
			\end{equation*}
			e, portanto, a			
			\begin{equation*}
				\boldsymbol{p_{\mathrm{im}}} \cdot \boldsymbol{f} + \boldsymbol{p_{\mathrm{ob}}} \cdot \boldsymbol{f} = \left( \boldsymbol{p_{\mathrm{im}}} + \boldsymbol{p_{\mathrm{ob}}} \right) \boldsymbol{f} = \boldsymbol{p_{\mathrm{ob}}} \cdot \boldsymbol{p_{\mathrm{im}}} \ .
			\end{equation*}
			
			Assim, notando não apenas que esse último resultado pode ser reescrito como			
			\begin{equation}
				\frac{\boldsymbol{p_{\mathrm{im}}} + \boldsymbol{p_{\mathrm{ob}}}}{\boldsymbol{p_{\mathrm{ob}}} \cdot \boldsymbol{p_{\mathrm{im}}}} = \frac{\boldsymbol{1}}{\boldsymbol{f}} \ , \label{quase}
			\end{equation}
			mas diante do fato que			
			\begin{equation}
				\frac{\boldsymbol{p_{\mathrm{im}}} + \boldsymbol{p_{\mathrm{ob}}}}{\boldsymbol{p_{\mathrm{ob}}} \cdot \boldsymbol{p_{\mathrm{im}}}} = \frac{\boldsymbol{p_{\mathrm{im}}}}{\boldsymbol{p_{\mathrm{ob}}} \cdot \boldsymbol{p_{\mathrm{im}}}} + \frac{\boldsymbol{p_{\mathrm{ob}}}}{\boldsymbol{p_{\mathrm{ob}}} \cdot \boldsymbol{p_{\mathrm{im}}}} = \frac{1}{\boldsymbol{p_{\mathrm{ob}}}} + \frac{1}{\boldsymbol{p_{\mathrm{im}}}} \ , \label{manipulation}
			\end{equation}
			é imediato concluir que a substituição de (\ref{manipulation}) na expressão (\ref{quase}) nos leva a
			\begin{equation*}
				\frac{\boldsymbol{1}}{\boldsymbol{p_{\mathrm{ob}}}} + \frac{\boldsymbol{1}}{\boldsymbol{p_{\mathrm{im}}}} = \frac{\boldsymbol{1}}{\boldsymbol{f}} \ ,
			\end{equation*}
			que é justamente a equação de Gauss cuja validade nós queríamos demonstrar.

\bigskip{\small \smallskip\noindent Updated: \today.}


\begin{thebibliography}{50}
		\bibitem{elon1} E. L. Lima: \emph{Qual o valor de $ 0^{0} $?} (Revista do Professor de Matemática nº 01, 1985).
		\bibitem{elon2} E. L. Lima: \emph{Novamente $ 0^{0} $} (Revista do Professor de Matemática 07, 1985).
		\bibitem{bancaria} P. Freire: \emph{Educação \textquotedblleft bancária\textquotedblright \hspace*{0.01cm} e educação libertadora}, \emph{Pedagogia do oprimido} (Rio de Janeiro, Paz e Terra, 1970) -- segunda edição, Capítulo II, p. 65-87.
		\bibitem{sutanto} S. H. Sutanto, P. C. Tjiang: \emph{J. Opt.} \textbf{13}, 105706 (2011).
		\bibitem{hecht} E. Hecht: \emph{Optics -- Fourth Edition} (Addison Wesley, San Francisco 2002).
		\bibitem{savelyev} I. V. Savelyev: \emph{Physics, A General Course -- Volume II} (Mir Publishers, Moscow 1989).
		\bibitem{moyses} H. M. Nussenzveig: \emph{Curso de Física Básica, Volume 4} (Edgard Blücher, Rio de Janeiro).
		\bibitem{piaget1} J. Piaget: \emph{Psicologia e Pedagogia, 4a edição} (Forense Universitária, Rio de Janeiro 1976).
		\bibitem{piaget2} J. Piaget: \emph{A formação do símbolo na criança} (Zahar, Rio de Janeiro 1978).
		\bibitem{piaget3} J. Piaget: \emph{Seis Estudos de Psicologia, 13a edição} (Forense Universitária, Rio de Janeiro 1985).
		\bibitem{photosforyou} https://pixabay.com/pt/photos/esférico-bola-universo-planeta-3140452/
		\bibitem{fireboltbyl} https://pixabay.com/pt/photos/espelho-esférico-espelho-edifício-1402472/
		\bibitem{pixel2013} https://pixabay.com/images/id-1576109/
		\bibitem{jjuni} https://pixabay.com/pt/illustrations/agora-medição-o-governante-1135297/
		\bibitem{openclip-1} https://pixabay.com/pt/vectors/xadrez-jogo-peças-peças-de-xadrez-145184/
		\bibitem{openclip-2} https://pixabay.com/pt/vectors/seta-red-brilhante-direito-próximo-145781/
		\bibitem{clker-free} https://pixabay.com/pt/vectors/vela-aceso-celebração-férias-cera-33287/
		\bibitem{user121799} https://tex.stackexchange.com/questions/481283/intersection-of-a-sphere-and-a-plane-knowing-equations
		\bibitem{crew} H. Crew: \emph{The Wave Theory of Light -- Memoirs of Huygens, Young and Fresnel} (American Book Company, 1900).
		\bibitem{newton-optica}  I. Newton: \emph{Óptica} (Edusp, São Paulo 2002).
		\bibitem{young-light-1} T. Young: \emph{Phil. Trans. R. Soc. Lond.} \textbf{92}, 12 (1802). 
		\bibitem{young-light-2} T. Young: \emph{Phil. Trans. R. Soc.} \textbf{94}, 1 (1804). 
		\bibitem{young} T. Young: \emph{Lectures on Natural Philosophy} (J. Johnson, London 1807).
		\bibitem{maxwell-wave} J. C. Maxwell: \emph{Phil. Trans. R. Soc.} \textbf{155}, 459 (1864). 
		\bibitem{hertz} M. de Abreu Faro: \emph{GFIS} \textbf{18} (2), 22 (1995).
		\bibitem{ribeiro} A. R. Ribeiro, L. Coelho, O. Bertolami, R. André: \emph{GFIS} \textbf{39} (1-2), 6 (2016). 
		\bibitem{salvetti} A. R. Salvetti: \emph{A História da Luz} (Editora Livraria da Física, São Paulo 2008).
		\bibitem{okuno} E. Okuno, I. L. Caldas, C. Chow: \emph{Física para Ciências Biológicas e Biomédicas} (Harper \& Row do Brasil, São Paulo 1982).
	\end{thebibliography}
\end{document}